\pdfoutput=1		

\documentclass[a4paper,12pt]{article}  
\usepackage[vmargin={20mm,20mm},hmargin={20mm,20mm}]{geometry}
\usepackage{commath}
\usepackage{graphicx}
\usepackage{enumitem}
\usepackage{tikz-cd}
\usepackage{subcaption}

\usepackage{amsmath}
\usepackage{amssymb}
\usepackage{amsthm}
\usepackage{bbm}
\usepackage{mathrsfs}

\let\emptyset\varnothing

\DeclareMathOperator{\Ad}{Ad}

\DeclareMathOperator{\Crit}{Crit}

\DeclareMathOperator{\Diff}{Diff}

\DeclareMathOperator{\grad}{grad}

\DeclareMathOperator{\Hess}{Hess}

\DeclareMathOperator{\id}{id}
\DeclareMathOperator{\im}{im}

\DeclareMathOperator{\ord}{ord}
\DeclareMathOperator{\pr}{pr}

\DeclareMathOperator{\supp}{supp}
\DeclareMathOperator{\Symp}{Symp}

\newcommand{\bld}[1]{\boldmath\textit{\textbf{#1}}\unboldmath}

\newtheoremstyle{main} 		             	 		
  	{}	                                     		
  	{}	                                    		
  	{\itshape}			                     		
  	{}        	                             		
  	{\boldmath\bfseries}   	                         		
  	{.}            	                        		
  	{ }           	                         		
  	{\thmname{#1}\thmnumber{ #2}\thmnote{ (#3)}}	
\theoremstyle{main}
\newtheorem{definition}{Definition}[section]
\newtheorem{axiom}{Axiom}[section]

\newtheorem{corollary}{Corollary}[section]
\newtheorem{theorem}{Theorem}[section]
\newtheorem{lemma}{Lemma}[section]
\newtheoremstyle{nonit} 		             	 		
  	{}	                                     		
  	{}	                                    		
  	{}			                     		
  	{}        	                             		
	{\boldmath\bfseries}   	                         		
  	{.}            	                        		
  	{ }           	                         		
  	{\thmname{#1}\thmnumber{ #2}\thmnote{ (#3)}}	
\theoremstyle{nonit}
\newtheorem{remark}{Remark}[section]
\newtheorem{example}{Example}[section]
\newtheoremstyle{ex} 		             	 		
  	{}	                                     		
  	{}	                                    		
  	{}			                     		
  	{}        	                             		
  	{\bfseries\boldmath}   	                         		
  	{.}            	                        		
  	{ }           	                         		
  	{\thmname{#1}\thmnumber{ #2}\thmnote{ (#3)}}	
\theoremstyle{ex}
\newtheorem{exercise}{Exercise}
\newtheorem{question}{Question}
\newtheorem*{solution}{Solution}

\usepackage[backend=bibtex, style=alphabetic]{biblatex}
\bibliography{bibliography}
\usepackage[bookmarksopen=true,bookmarksnumbered=true]{hyperref}
\hypersetup{
  colorlinks   = true, 
  urlcolor     = blue, 
  linkcolor    = blue, 
  citecolor    = blue 
}

\usepackage{tcolorbox}
\tcbuselibrary{theorems}
\newtcbtheorem
  []
  {summary}
  {Summary}
  {%
    colback=blue!5,
    colframe=blue!35!black,
    fonttitle=\bfseries,
  }
  {sum}

\begin{document}

\title{Lectures on Classical Mechanics\\
\emph{A Didactical Approach to Higher Mathematics}
}


\author{Yannis B\"ahni
}

\maketitle

\begin{abstract}
	The aim of this paper is twofold: First, we give a formal introduction to the basics of the mathematical framework of classical mechanics. Along the way, we prove a Hamiltonian and a Lagrangian version of Noether's Theorem, an important result concerning continuous symmetries of physical systems. At the end, we prove a new statement about orbit cylinders on homotopies of stable regular energy hypersurfaces. The main question we answer is \emph{what is the dynamical meaning of stability?} The second aim is to provide a didactical framework for introducing advanced mathematics, which can also be of used in other topics. A broad range of established methods belonging to the realm of cognitive activation as well as formative assessment techniques are used.	
\end{abstract}

\tableofcontents

\newpage
\section{Introduction}
\subsection{Symplectic Geometry}
The modern language of classical mechanics is provided by symplectic geometry. For an introduction to symplectic geometry see \cite{silva:sg:2008} and for a more sophisticated treatment \cite{mcduffsalamon:st:2017}. For an introduction to Hamiltonian dynamics see \cite{abrahammarsden:cm:1978} as well as \cite{hoferzehnder:hd:1994} for a view towards symplectic invariants. Moreover, the books \cite{arnold:cm:2006} and \cite{cushman:cm:2015} provide a much more in depth introduction to the mathematical treatment of classical mechanics. As a short introduction we recommend \cite[Chapter~1]{takhtajan:qm:2008}. The only viewpoint of symplectic geometry we want to elaborate on, is its drastic contrast to Riemannian geometry. The following exercise is inspired by \cite{gotay:sg:1992}. 

\begin{figure}[h!tb]
     \centering
     \begin{subfigure}[b]{0.47\textwidth}
         \centering
         \includegraphics[width=\textwidth]{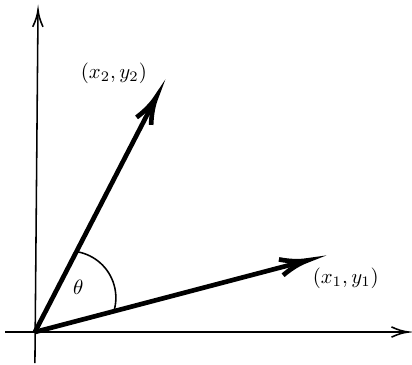}
     \end{subfigure}
     \hfill
     \begin{subfigure}[b]{0.47\textwidth}
         \centering
         \includegraphics[width=\textwidth]{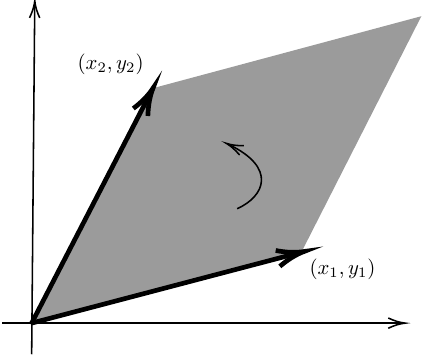}
     \end{subfigure}
	 \caption{}
	 \label{fig:anglearea}
\end{figure}

\begin{exercise}
	\label{ex:introduction}
	On the real vector space $\mathbb{R}^2$ we consider the two functions
	\begin{equation*}
		\langle \cdot, \cdot \rangle \colon \mathbb{R}^2 \times \mathbb{R}^2 \to \mathbb{R}, \qquad \langle (x_1,y_1), (x_2,y_2) \rangle := x_1x_2 + y_1 y_2, 
	\end{equation*}
	\noindent and
	\begin{equation*}
		\omega\colon \mathbb{R}^2 \times \mathbb{R}^2 \to \mathbb{R}, \qquad \omega((x_1,y_1),(x_2,y_2)) := x_1 y_2 - x_2 y_1. 
	\end{equation*}

	\begin{enumerate}[label=\textup{(\alph*)}]
		\item What is the geometrical meaning of the two functions $\langle \cdot,\cdot \rangle$ and $\omega$?

		\item What are common properties of the two functions $\langle\cdot,\cdot \rangle$ and $\omega$? Which properties are different? It might help to consider the induced maps
			\begin{equation*}
				\langle \cdot, \cdot \rangle \colon \mathbb{R}^2 \vee \mathbb{R}^2 \to \mathbb{R}, \qquad \langle (x_1,y_1) \vee (x_2,y_2) \rangle = x_1x_2 + y_1 y_2, 
			\end{equation*}
			\noindent and
			\begin{equation*}
				\omega \colon \mathbb{R}^2 \wedge \mathbb{R}^2 \to \mathbb{R}, \qquad \omega((x_1,y_1) \wedge (x_2,y_2)) = x_1 y_2 - x_2 y_1.
			\end{equation*}

			\item Consider the induced functions
			\begin{equation*}
				\Phi \colon \mathbb{R}^2 \to (\mathbb{R}^2)^*, \qquad \Phi(x_1,y_1)(x_2,y_2) := \langle (x_1,y_1), (x_2,y_2) \rangle,
			\end{equation*}
			\noindent and
			\begin{equation*}
				\Psi \colon \mathbb{R}^2 \to (\mathbb{R}^2)^*, \qquad \Psi(x_1,y_1)(x_2,y_2) := \omega((x_1,y_1),(x_2,y_2)).
			\end{equation*}
			Which properties do these two functions share? Are there any different properties? \emph{Hint}: In particular, check if $\Phi$ and $\Psi$ are injective or surjective functions.
		\item What is the geometrical meaning of the properties of $\Phi$ and $\Psi$?
	\end{enumerate}
\end{exercise}

\begin{solution}
	\mbox{}
	\begin{enumerate}[label=\textup{(\alph*)}]
		\item The function $\langle \cdot, \cdot \rangle$ is the standard inner product on $\mathbb{R}^2$. We have that
			\begin{equation*}
				\langle a,b \rangle = \norm[0]{a}\norm[0]{b} \cos \theta \qquad \forall a,b \in \mathbb{R}^2,
			\end{equation*}
			\noindent where $\theta$ is the angle between the vectors $a$ and $b$. This is the left illustration in Figure \ref{fig:anglearea}. For $\omega$ we observe
            \begin{equation*}
				\omega(a,b) = \det(a\vert b) \qquad \forall a,b \in \mathbb{R}^2.
			\end{equation*}
			Thus $\omega(a,b)$ is the signed area of the parallelogram spanned by the two vectors $a$ and $b$. This is the illustration on the right in Figure \ref{fig:anglearea}.
		\item The two functions $\langle \cdot,\cdot \rangle$ and $\omega$ are both bilinear forms and thus induce unique linear functions
			\begin{equation*}
				\langle \cdot, \cdot \rangle \colon \mathbb{R}^2 \otimes \mathbb{R}^2 \to \mathbb{R}, \qquad \langle (x_1,y_1) \otimes (x_2,y_2) \rangle = x_1x_2 + y_1 y_2, 
			\end{equation*}
			\noindent and
			\begin{equation*}
				\omega \colon \mathbb{R}^2 \otimes \mathbb{R}^2 \to \mathbb{R}, \qquad \omega((x_1,y_1) \otimes (x_2,y_2)) = x_1 y_2 - x_2 y_1.
			\end{equation*}
			\noindent The bilinear form $\langle \cdot,\cdot \rangle$ is symmetric and so induces a unique linear map
						\begin{equation*}
				\langle \cdot, \cdot \rangle \colon \mathbb{R}^2 \vee \mathbb{R}^2 \to \mathbb{R}, \qquad \langle (x_1,y_1) \vee (x_2,y_2) \rangle = x_1x_2 + y_1 y_2, 
			\end{equation*}
			\noindent where $\mathbb{R}^2 \vee \mathbb{R}^2$ denotes the symmetric tensor product of $\mathbb{R}^2$. The bilinear form $\omega$ is antisymmetric and so induces a unique linear map
			\begin{equation*}
				\omega \colon \mathbb{R}^2 \wedge \mathbb{R}^2 \to \mathbb{R}, \qquad \omega((x_1,y_1) \wedge (x_2,y_2)) = x_1 y_2 - x_2 y_1
			\end{equation*}
			\noindent where $\mathbb{R}^2 \wedge \mathbb{R}^2$ denotes the alternating tensor product of $\mathbb{R}^2$. The bilinear form $\langle \cdot, \cdot \rangle$ is positive definite and $\omega$ has the property that $\omega(a,a) = 0$ for all $a \in \mathbb{R}^2$.
		\item The induced maps $\Phi$ and $\Psi$ are well defined maps to the dual space $(\mathbb{R}^2)^* \cong \mathbb{R}^2$ because of part b. Both maps are bijective. Indeed, if $\Phi(x_1,y_1) = \Phi(x_2',y_2')$, then 
			\begin{equation*}
				\Phi(x_1,y_1)(x_2,y_2) = \Phi(x_1',y_1')(x_2,y_2) \qquad \forall (x_2,y_2) \in \mathbb{R}^2. 
			\end{equation*}
			In particular, it holds that
			\begin{equation*}
				x_1 = \Phi(x_1,y_1)(1,0) = \Phi(x_1',y_1')(1,0) = x_1' 
			\end{equation*}
			\noindent and
			\begin{equation*}
				y_1 = \Phi(x_1,y_1)(0,1) = \Phi(x_1',y_1')(0,1) = y_1'. 
			\end{equation*}
			Hence $\Phi$ is injective. An injective linear map between finite-dimensional vector spaces of the same dimension is automatically a bijection. Analogously, one shows that $\Psi$ is bijective.
		\item Two nonzero vectors $a, b \in \mathbb{R}^2$ are orthogonal, if and only if $\langle a, b \rangle = 0$ holds. The bijectivity of $\Phi$ is equivalent to the statement, that a vector $a \in \mathbb{R}^2$ cannot be orthogonal to every vector $b \in \mathbb{R}^2$, unless $a = 0$. The bijectivity of $\Psi$ geometrically means that only the degenerate parallelogram has no area, meaning that $a$ and $b$ are linearly dependent. 
	\end{enumerate}
\end{solution}

Mathematics is known for its high degree of abstraction. By \cite[p.~37]{stern:lp:2021}, many students fail at this high level of abstractness already at the introduction of Algebra. Why are such abstract concepts so hard to grasp? In contrast to concepts like ``creature'', mathematical concepts most often cannot be reduced to concrete, sensually perceptible events \cite[p.~39]{stern:lp:2021}. Teaching Mathematics, and more broadly, subjects belonging to the MINT regime, is hard. A didactical step towards better teaching methods on the high-school level was undertaken in \cite{stern:lp:2022}. The aim of this work is to show, that these methods can also be implemented on a much higher level of mathematics, namely research level. As is this field is far to vast, we focus on the authors very recent paper \cite{baehni:rfh:2023}. The aim of this paper was to construct a generalisation of Rabinowitz--Floer homology. See the excellent survey article \cite{albersfrauenfelder:rfh:2012} for more details. Our aim is twofold:
\begin{itemize}
	\item Performing a feasibility analysis whether results of teaching and learning research can be efficiently applied in such a complex scenario.
	\item Laying the basis for the book project \emph{Lectures on Twisted Rabinowitz--Floer Homology: A Didactical Approach to Higher Mathematics}. To the authors knowledge, there is currently no book treating Rabinowitz--Floer homology except the excellent foundational papers as well as the book project \cite{weber:orbits:2018}.
\end{itemize}

Since prior knowledge is among the best predictors for learning \cite[Section~2]{ralph:knowledge:2023}, we carry out a carefully designed formative assessment to reactivate student knowledge of differential topology and symplectic geometry as well as tackling misconceptions. This is an implementation of the constructivism approach to learning, where new knowledge builds upon prior knowledge by replacing naive models with expert models.

The Lie derivative generalises the directional derivative of a function to arbitrary tensor fields on a smooth manifold. Let $M$ be a smooth manifold and $A \in \mathcal{T}^{k,l}(M)$ a tensor field of type $(k,l)$. Then for any $X \in \mathfrak{X}(M)$ we define the \bld{Lie derivative of $A$ with respect to $X$} to be the tensor field $L_X A \in \mathcal{T}^{k,l}(M)$ given by
	\begin{equation*}
		L_X A := \frac{d}{dt}\bigg\vert_{t = 0} \phi^*_t A,
	\end{equation*}
    \noindent where $\phi$ denotes the smooth flow of $X$. Equivalently, the Lie derivative is the unique tensor derivation such that
	\begin{equation*}
		L_X f = Xf \qquad \text{and} \qquad L_XY = [X,Y]
	\end{equation*}
	\noindent for all $f \in C^\infty(M)$ and $X,Y \in \mathfrak{X}(M)$, where $[X,Y] = XY - YX$ denotes the Lie bracket on the Lie algebra of vector fields $\mathfrak{X}(M)$.

    The exterior differential generalises the differential of a function to arbitrary differential forms. Let $\omega \in \Omega^k(M)$ be a differential $k$-form written locally as $\omega = \sum_I \omega_I dx_I$, where the sum runs over all increasing multiindices $I$. Then we define the \bld{exterior differential of $\omega$} to be the differential $k + 1$-form $d\omega \in \Omega^{k + 1}(M)$ given locally by
	\begin{equation*}
		d\omega = \sum_I d\omega_I \wedge dx_I.
	\end{equation*}
    \noindent Equivalently, the exterior differential is the unique graded derivation of degree $1$ such that 
    \begin{equation*}
        df(X) = Xf = L_X f \qquad \text{and} \qquad d \circ d = 0
    \end{equation*}
    \noindent for all $f \in C^\infty(M)$ and $X \in \mathfrak{X}(M)$.

    \begin{lemma}[{Cartan's Magic Formula, \cite[Theorem~14.35]{lee:dt:2012}}]
        \label{lem:cartan}
        Let $M$ be a smooth manifold. Then
        \begin{equation*}
            L_X = i_X \circ d + d \circ i_X
        \end{equation*}
        \noindent holds for all $X \in \mathfrak{X}(M)$, where
        \begin{equation*}
            i_X \colon \Omega^{* + 1}(M) \to \Omega^*(M), \qquad i_X\omega := \omega(X,\cdot,\hdots,\cdot)
        \end{equation*}
        \noindent denotes interior multiplication.     
    \end{lemma}

    \begin{lemma}[{Fisherman's Formula, \cite[Proposition~22.14]{lee:dt:2012}}]
        \label{lem:fisherman}
        Let $M$ be a smooth manifold. Then 
        \begin{equation*}
            \frac{d}{dt} \phi_t^*\omega = \phi_t^*L_{X_t}\omega
        \end{equation*}
        \noindent holds for all $\omega \in \Omega(M)$ and time-dependent vector fields $X$.
    \end{lemma}

\begin{question}
        Let $M$ be a smooth manifold. Which statements are correct on $\Omega(M)$?
    \begin{itemize}
        \item[$\square$] $L_X (\omega \wedge \eta) = L_X\omega \wedge \eta + \omega \wedge L_X\eta$ for all $\omega \in \Omega^k(M)$, $\eta \in \Omega^l(M)$ and $X \in \mathfrak{X}(M)$.
        \item[$\square$] $L_X (\omega \wedge \eta) = L_X\omega \wedge \eta + (-1)^{kl}\omega \wedge L_X\eta$ for all $\omega \in \Omega^k(M)$, $\eta \in \Omega^l(M)$ and $X \in \mathfrak{X}(M)$.
        \item[$\square$] $d (\omega \wedge \eta) = d\omega \wedge \eta + \omega \wedge d\eta$ for all $\omega \in \Omega^k(M)$ and $\eta \in \Omega^l(M)$.
        \item[$\square$] $d(\omega \wedge \eta) = d\omega \wedge \eta + (-1)^{kl}\omega \wedge d\eta$ for all $\omega \in \Omega^k(M)$ and $\eta \in \Omega^l(M)$.
    \end{itemize}
\end{question}

\begin{solution}
    Correct is the first and the fourth item.
\end{solution}

\begin{question}
    Let $M$ be a smooth manifold. Which statements are correct on $\Omega(M)$?
    \begin{itemize}
        \item[$\square$] $i_{[X,Y]} = L_Yi_X - i_XL_Y$ for all $X,Y\in \mathfrak{X}(M)$.
        \item[$\square$] $i_{[X,Y]} = L_Xi_Y - i_YL_X$ for all $X,Y\in \mathfrak{X}(M)$.
        \item[$\square$] $i_X (\omega \wedge \eta) = i_X\omega \wedge \eta + \omega \wedge i_X\eta$ for all $\omega \in \Omega^k(M)$, $\eta \in \Omega^l(M)$ and $X \in \mathfrak{X}(M)$.
        \item[$\square$] $i_X(\omega \wedge \eta) = i_X\omega \wedge \eta + (-1)^{kl}\omega \wedge i_X\eta$ for all $\omega \in \Omega^k(M)$, $\eta \in \Omega^l(M)$ and $X \in \mathfrak{X}(M)$.
    \end{itemize}
\end{question}

\begin{solution}
    Correct is the second and the fourth item.
\end{solution}

\begin{question}
    Let $M$ be a smooth manifold. A differential form $\omega \in \Omega^2(M)$ on $M$ is called a \bld{symplectic form}, if
    \begin{itemize}
        \item[$\square$] $\omega$ is nondegenerate, meaning that
        \begin{equation*}
            \Phi_X \colon \mathfrak{X}(M) \to \Omega^1(M), \qquad \Phi_X(Y) := \omega(X,Y)
        \end{equation*}
        \noindent is an isomorphism for all vector fields $X \in \mathfrak{X}(M)$.
        \item[$\square$] $\omega$ is nondegenerate and exact, meaning that there exists a differential form $\lambda \in \Omega^1(M)$ with $\omega = d\lambda$ for the exterior differential $d$.
        \item[$\square$] $\omega$ is nondegenerate and closed, meaning $d\omega = 0$.
        \item[$\square$] the dimension $\dim M$ of $M$ is even, say $\dim M = 2n$, and $\omega^n \in \Omega^{2n}(M)$ is a volume form, meaning that $\omega^n$ is nowhere vanishing.
    \end{itemize}
\end{question}

\begin{solution}
    Correct is the third and fourth item.
\end{solution}

\begin{question}
    \label{que:symplectic}
    Let $(M^{2n},\omega)$ be a symplectic manifold. Which statements are correct?
    \begin{itemize}
        \item[$\square$] $H^2_{\operatorname{dR}}(M) \neq 0$, where $H^2_{\operatorname{dR}}(M)$ denotes the second de Rham cohomology group of $M$.
        \item[$\square$] If $M$ is compact, then $[\omega] \neq 0$ in $H^2_{\operatorname{dR}}(M)$, meaning that $\omega$ is not exact.
        \item[$\square$] No even sphere $\mathbb{S}^{2n - 1}$ does admit a symplectic form for $n \geq 2$.
        \item[$\square$] The cotangent bundle of any smooth manifold is not canonically orientable.
    \end{itemize}
\end{question}

\begin{solution}
    Correct is the second and the third item.
\end{solution}

\begin{question}
    Let $\omega \in \Omega^2(M)$ be nondegenerate on a smooth manifold $M^{2n}$. Which statements are correct?
    \begin{itemize}
        \item[$\square$] For every point $(x,y) \in M$ there does exist a chart $\del[1]{U,(x_i,y_i)}$ about $(x,y)$ with
        \begin{equation*}
            \omega\vert_U = \sum_{i = 1}^n dy_i \wedge dx_i.
        \end{equation*}
        Such coordinates are called \bld{Darboux coordinates}.
        \item[$\square$] Like Riemannian manifolds, symplectic manifolds do carry local invariants. For example, a local invariant in Riemannian geometry is curvature.
        \item[$\square$] Unlike Riemannian manifolds, symplectic manifolds do not carry any local invariants and are thus locally indistinguishable.
        \item[$\square$] The only invariants worth studying in symplectic geometry are global invariants like symplectic capacities or Floer homology groups.
    \end{itemize}
\end{question}

\begin{solution}
    Correct is the third and the fourth item.
\end{solution}

\subsection{Periodic Orbits}
According to \cite[p.~169]{albersfrauenfelder:leafwise:2012}, periodic orbits are the fundamental building blocks for Hamiltonian systems. For example, the existence of closed orbits is important in the study of celestial mechanics. Indeed, consider the restricted three body problem, that is, a massless satellite moving in the gravitational field of two planets. In order to search for extraterrestrial life, one needs to know if the satellite can orbit around one of them. However, finding periodic orbits in whatever sense is no easy task in general. For getting a feeling, we first perform some observations in a very geometrical setting. 

\begin{exercise}
	\label{ex:orbits_sphere}
	Consider the odd-dimensional sphere 
	\begin{equation*}
		\mathbb{S}^{2n - 1} := \cbr[4]{(z_1,\dots,z_n) \in \mathbb{C}^n \colon \sum_{j = 1}^n \abs[0]{z_j}^2 = 1}.
	\end{equation*}
	Compute the flow of the vector field 
	\begin{equation*}
		2\sum_{j = 1}^n \del[4]{y_j\frac{\partial}{\partial x_j} - x_j \frac{\partial}{\partial y_j}} \in \mathfrak{X}(\mathbb{R}^{2n})
	\end{equation*}
	\noindent on $\mathbb{S}^{2n - 1}$, where $z_j := x_j + iy_j$, and conclude that this flow is periodic.
\end{exercise}

\begin{solution}
    Writing the vector field in the complex coordinates $z_j$ yields the vector field
    \begin{equation*}
        2i \sum_{j = 1}^n \del[4]{\bar{z}_j\frac{\partial}{\partial \bar{z}_j} - z_j  \frac{\partial}{\partial z_j}}
    \end{equation*}
    \noindent with
    \begin{equation*}
        \frac{\partial}{\partial z_j} = \frac{1}{2}\del[3]{\frac{\partial}{\partial x_j} - i \frac{\partial}{\partial y_j}} \qquad \text{and} \qquad \frac{\partial}{\partial \bar{z}_j} = \frac{1}{2}\del[3]{\frac{\partial}{\partial x_j} + i \frac{\partial}{\partial y_j}}.
    \end{equation*}
    Thus the periodic flow is given by
    \begin{equation*}
        \mathbb{R} \times \mathbb{S}^{2n - 1} \to \mathbb{S}^{2n - 1}, \qquad (t,z) \mapsto e^{-2it}z.
    \end{equation*}
\end{solution}

\begin{exercise}
	\label{ex:orbits_ellipsoid}
	For real numbers $a_1,\dots,a_n > 0$, consider the ellipsoids
	\begin{equation*}
		E(a_1,\dots,a_n) := \cbr[4]{(z_1,\dots,z_n) \in \mathbb{C}^n \colon \sum_{j = 1}^n \frac{\pi\abs[0]{z_j}^2}{a_j} = 1}.
	\end{equation*}
	Compute the flow of the vector field
	\begin{equation*}
		2\sum_{j = 1}^n \frac{\pi}{a_j}\del[4]{y_j\frac{\partial}{\partial x_j} - x_j \frac{\partial}{\partial y_j}} \in \mathfrak{X}(\mathbb{R}^{2n})
	\end{equation*}
	\noindent on $E(a_1,\dots,a_n)$ and determine the closed orbits.
\end{exercise}

\begin{solution}
    In complex coordinates $z_j = x_j + iy_j$ the vector field is given by
    \begin{equation*}
        2i \sum_{j = 1}^n \frac{\pi}{a_j}\del[4]{\bar{z}_j\frac{\partial}{\partial \bar{z}_j} - z_j  \frac{\partial}{\partial z_j}}.
    \end{equation*}
    Thus the flow is
    \begin{equation*}
        \mathbb{R} \times E(a_1,\dots,a_n) \to E(a_1,\dots,a_n), \qquad (t,z) \mapsto \del[1]{e^{-2i \pi t/a_1}z_1,\dots,e^{-2\pi i t/a_n}z_n}.
    \end{equation*}
    Now by \cite[p.~26]{schlenk:embedding:2005}, the periodic orbits depend on the numbers $a_j$. If the numbers $a_j$ are linearly independent over $\mathbb{Z}$, the only periodic orbits are of the form
    \begin{equation*}
        \mathbb{R} \to E(a_1,\dots,a_n), \qquad t \mapsto \del[1]{0,\dots, e^{-2\pi i t/a_j}z_j, \dots, 0}.
    \end{equation*}
    The general case is slightly more difficult.
\end{solution}

\begin{exercise}
	\label{ex:orbits_star-shaped}
	Let $\Sigma \subseteq \mathbb{C}^n$ be a compact and connected smooth star-shaped hypersurface with respect to the origin $0 \in \mathbb{C}^n$. Define a function
	\begin{equation*}
		h_\Sigma \colon \mathbb{C}^n \setminus \{0\} \to \mathbb{R}, \qquad h_\Sigma(z) := \frac{1}{r(z)^2},
	\end{equation*}
	\noindent where
	\begin{equation*}
		r \colon \mathbb{C}^n \setminus \{0\} \to \intoo[0]{0,+\infty} 
	\end{equation*}
	\noindent is the unique smooth function such that $r(z)z \in \Sigma$. The induced vector field is given by
	\begin{equation}
		\label{eq:Reeb_star-shaped}
		X_\Sigma = \sum_{j = 1}^n \del[4]{\frac{\partial h_\Sigma}{\partial y_j}\frac{\partial}{\partial x_j} - \frac{\partial h_\Sigma}{\partial x_j}\frac{\partial}{\partial y_j}} \in \mathfrak{X}(\mathbb{R}^{2n}).
	\end{equation}
	\begin{enumerate}[label=\textup{\alph*.}]
		\item Show, that the vector field $X_{\mathbb{S}^{2n - 1}}$ extends to the vector field given in Exercise \ref{ex:orbits_sphere}.
		\item Show, that the vector field $X_{E(a_1,\dots,a_n)}$ extends to the vector field given in Exercise \ref{ex:orbits_ellipsoid}.
		\item Explain, why every such hypersurface can be written as
			\begin{equation*}
				\Sigma_f := \cbr[0]{f(z)z : z \in \mathbb{S}^{2n - 1}}
			\end{equation*}
			\noindent for a suitable positive function $f \in C^\infty(\mathbb{S}^{2n - 1})$.
		\item Given $f \in C^\infty(\mathbb{S}^{2n - 1})$ positive, compute the corresponding vector field $X_{\Sigma_f}$. Does $X_{\Sigma_f}$ admit a closed orbit on $\Sigma_f$? Can you prove your guess? 
	\end{enumerate}
\end{exercise}

\begin{figure}[h!tb]
	\centering
	\includegraphics[width=\textwidth]{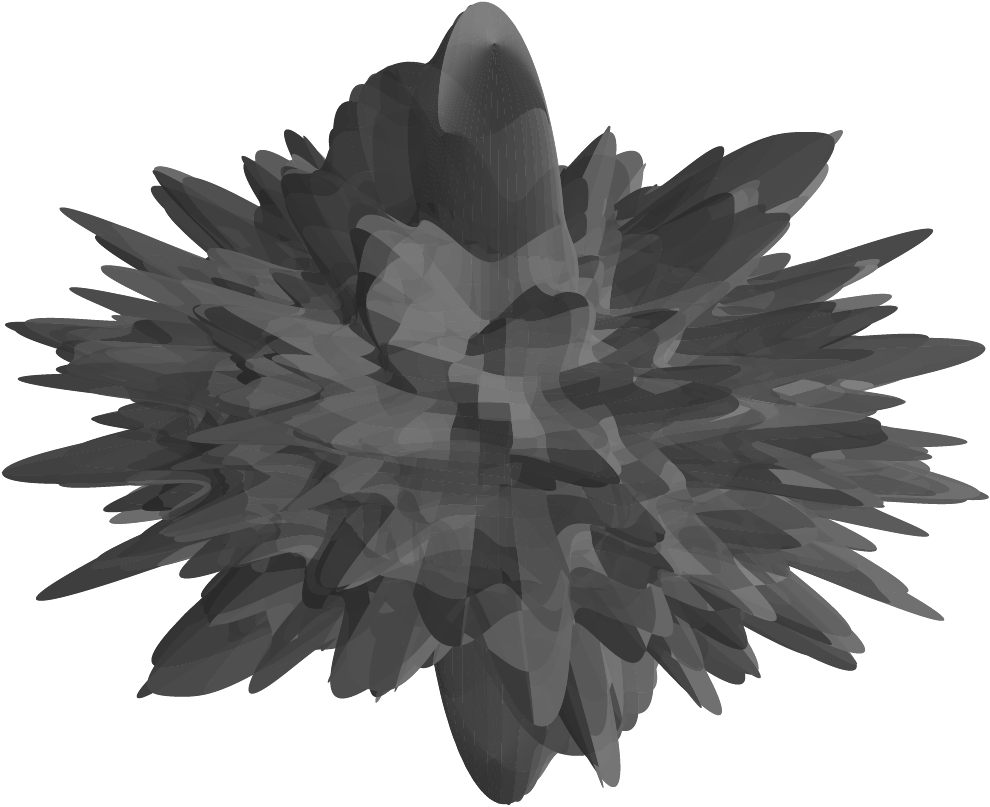}
	\caption{A $\mathbb{Z}_2$-symmetric star-shaped hypersurface $\Sigma_f$}
	\label{fig:spiked_sphere}
\end{figure}

\begin{solution}
    \mbox{}
    \begin{enumerate}[label=\alph*.]
        \item In this case, we have that
            \begin{equation*}
                h_{\mathbb{S}^{2n - 1}} \colon \mathbb{C}^n \setminus \{0\} \to \mathbb{R}, \qquad h_{\mathbb{S}^{2n - 1}}(z) = \norm[0]{z}^2
            \end{equation*}
            \noindent and thus
            \begin{equation*}
                X_{\mathbb{S}^{2n - 1}} = 2\sum_{j = 1}^n \del[4]{y_j\frac{\partial}{\partial x_j} - x_j \frac{\partial}{\partial y_j}}.
            \end{equation*}
        \item In this case, we have that
            \begin{equation*}
                h_{E(a_1,\dots,a_n)} \colon \mathbb{C}^n \setminus \{0\} \to \mathbb{R}, \qquad h_{E(a_1,\dots,a_n)}(z) = \sum_{j = 1}^n \frac{\pi \abs[0]{z_j}^2}{a_j}
            \end{equation*}
            \noindent and thus
            \begin{equation*}
                X_{\mathbb{S}^{2n - 1}} = 2\sum_{j = 1}^n \frac{\pi}{a_j}\del[4]{y_j\frac{\partial}{\partial x_j} - x_j \frac{\partial}{\partial y_j}}.
            \end{equation*}
        \item We immediately see that $\Sigma = \Sigma_r$.
        \item In this case, we have that
            \begin{equation*}
                h_{\Sigma_f} \colon \mathbb{C}^n \setminus \{0\} \to \mathbb{R}, \qquad h_{\Sigma_f}(z) = \norm[0]{z}^2 f\del[3]{\frac{z}{\norm[0]{z}}}^{-2}
            \end{equation*}
            \noindent and thus $X_{\Sigma_f}$ is fairly ugly to compute. One concludes from part a and b, that there exists a periodic orbit for this vector field on the sphere $\mathbb{S}^{2n - 1}$ and the ellipsoids $E(a_1,\dots,a_n)$. However, for general $f$, it is not very clear why there should always exist a periodic orbit by considering the wild star-shaped hypersurface in Figure \ref{fig:spiked_sphere}.
    \end{enumerate}
\end{solution}
\newpage
\section{The Hamiltonian Formalism}
The Hamiltonian description of classical mechanics admits a natural formulation in the language of symplectic geometry. It is one of two modern formulations of the laws of classical physics and translates also to the formulation of quantum mechanics. Indeed, the time-dependent \emph{Schr\"odinger equation} is given by
\begin{equation*}
    i\hbar \frac{d\psi}{dt} = H\psi,
\end{equation*}
\noindent where $H$ is a self-adjoint operator on an infinite-dimensional separable complex Hilbert space.

\subsection{Classical Observables}

\begin{definition}[Hamiltonian System]
	A \bld{Hamiltonian system}\index{Hamiltonian!system} is a symplectic manifold $(M,\omega)$, called the \bld{phase space}, together with a smooth function $H \in C^\infty(M)$, called a \bld{Hamiltonian function}. We write $(M,\omega,H)$ for a Hamiltonian system.
\end{definition}

\begin{example}[Magnetic Hamiltonian System]
	Let $(M,m)$ be a pseudo-Riemannian manifold and denote by $\pi_{T^*M} \colon T^*M \to M$ its cotangent bundle. For a smooth potential function $V \in C^\infty(M)$ define a Hamiltonian function $H \in C^\infty(T^*M)$ by 
	\begin{equation}
		\label{eq:mechanical_Hamiltonian}
		H(q,p) := \frac{1}{2}\norm[0]{p}^2_{m^*} + V(q).
	\end{equation}
	For $\sigma \in \Omega^2(M)$ closed, the form $\omega_\sigma := \lambda + \pi^*_{T^*M}\sigma$ is a symplectic form on $T^*M$, where $(q,p)$ denote the standard coordinates on the cotangent bundle and $\lambda \in \Omega^1(T^*M)$ is the \bld{canonical Liouville form} defined by
    \begin{equation}
        \label{eq:liouville}
        \lambda_{(q,p)}(v) := p\del[1]{D\pi_{T^*M}(v)} \qquad \forall v \in T_{(q,p)}T^*M.
    \end{equation}
    Locally, we have that $\lambda = dp \wedge dq$. The symplectic manifold $(T^*M,\omega_\sigma)$ is called a \bld{magnetic cotangent bundle} and the Hamiltonian system $(T^*M,\omega_\sigma,H)$ is called a \bld{magnetic Hamiltonian system}. If $\sigma = 0$, the system is called a \bld{mechanical Hamiltonian system}.
\end{example}

\begin{definition}[Hamiltonian Vector Field]
	Let $(M,\omega,H)$ be a Hamiltonian system. The \bld{Hamiltonian vector field}\index{Hamiltonian!vector field} is defined to be the vector field $X_H \in \mathfrak{X}(M)$ given implicitly by
	\begin{equation*}
		i_{X_H}\omega = -dH.
	\end{equation*}	
\end{definition}

\begin{example}[Magnetic Hamiltonian Systems]
 	Let $(T^*M^n,\omega_\sigma,H)$ be a magnetic Hamiltonian system. Then the Hamiltonian vector field $X_H$ is given by 
	\begin{equation}
	\label{eq:magnetic_Hamiltonian_vf}
		X_H(q,p) = \sum_{i = 1}^n\frac{\partial H}{\partial p_i}\frac{\partial}{\partial q_i} + \sum_{i = 1}^n\del[4]{\sum_{j = 1}^n\sigma_{ij}(q)\frac{\partial H}{\partial p_j} - \frac{\partial H}{\partial q_i}}\frac{\partial}{\partial p_i},
	\end{equation}
	\noindent where
	\begin{equation*}
		\sigma = \frac{1}{2} \sum_{i,j = 1}^n\sigma_{ij}(q) dq_i \wedge dq_j, \qquad \sigma_{ji} = -\sigma_{ij}.
	\end{equation*}
\end{example}

\begin{exercise}
    Let $(M^{2n},\omega)$ be a symplectic manifold. Write down the local formula for the Hamiltonian vector field in Darboux coordinates $(q_j,p_j)$. Which equations do integral curves of the Hamiltonian vector field solve?
\end{exercise}

\begin{solution}
    In Darboux coordinates we have
    \begin{equation*}
        \omega = \sum_{j = 1}^n dp_j \wedge dq_j.
    \end{equation*}
    Thus the Hamiltonian vector field $X_H$ is locally given by
    \begin{equation*}
        X_H = \sum_{j = 1}^n\del[3]{\frac{\partial H}{\partial p_j}\frac{\partial}{\partial q_j} - \frac{\partial H}{\partial q_j}\frac{\partial}{\partial p_j}}.
    \end{equation*}
    Consequently, integral curves of $X_H$ satisfy \bld{Hamilton's equations}
    \begin{equation}
        \label{eq:hamiltonian_equations}
        \dot{q} = \frac{\partial H}{\partial p} \qquad \text{and} \qquad \dot{p} = -\frac{\partial H}{\partial q}.
    \end{equation}
\end{solution}

\begin{lemma}[{Jacobi, \cite[Theorem~3.3.19]{abrahammarsden:cm:1978}}]
	\label{lem:Jacobi}
	Let $(M,\omega,H)$ be a Hamiltonian system and let $\varphi \in \Symp(M,\omega)$ be a symplectomorphism. Then
	\begin{equation*}
		\varphi^*X_H = X_{\varphi^*H}.
	\end{equation*}
\end{lemma}

\begin{proof}
	We compute
	\begin{align*}
		i_{X_{\varphi^*H}}\omega = -d\varphi^*H = -\varphi^*dH = \varphi^*(i_{X_H}\omega) = i_{\varphi^*X_H}(\varphi^*\omega) = i_{\varphi^*X_H}\omega.
	\end{align*}
	Thus we conclude by the uniqueness of the Hamiltonian vector field.	
\end{proof}

\begin{lemma}
	\label{lem:conjugated_flow}
	Let $(M,\omega,H)$ be a Hamiltonian system and let $\varphi \in \Symp(M,\omega)$ be a symplectomorphism. Then 
	\begin{equation*}
		\phi^{X_{\varphi^*H}}_t = \varphi^{-1} \circ \phi_t^{X_H} \circ \varphi, 
	\end{equation*}
	\noindent whenever either side is defined, where $\phi$ denotes the smooth flow of a vector field.
\end{lemma}

\begin{proof}
	Using Lemma \ref{lem:Jacobi} we compute
	\begin{align*}
		\frac{d}{dt}\varphi^{-1} \circ \phi_t^{X_H} \circ \varphi &= D\varphi^{-1} \circ \frac{d}{dt}\phi^{X_H}_t \circ \varphi\\
		&= D\varphi^{-1} \circ X_H \circ \phi_t^{X_H} \circ \varphi\\
		&= D\varphi^{-1} \circ X_H \circ \varphi \circ \varphi^{-1} \circ \phi_t^{X_H} \circ \varphi\\
		&= \varphi^*X_H \circ \varphi^{-1} \circ \phi_t^{X_H} \circ \varphi\\
		&= X_{\varphi^*H} \circ \varphi^{-1} \circ \phi_t^{X_H} \circ \varphi,
	\end{align*}
	\noindent and the result follows by the uniqueness of integral curves.
\end{proof}

\begin{definition}[{Algebra of Classical Observables, \cite[p.~46]{takhtajan:qm:2008}}]
	Let $(M,\omega)$ be a symplectic manifold. The commutative real algebra $C^\infty(M)$ of smooth functions on $M$ is called the \bld{algebra of classical observables}.
\end{definition}

\begin{remark}
	In quantum mechanics, the algebra of observables is not commutative. This gives rise to very interesting phenomenons. For more details see \cite[p.~66]{takhtajan:qm:2008}.
\end{remark}

\begin{definition}[{Poisson Bracket, \cite[p.~578]{lee:dt:2012}}]
	Let $(M,\omega)$ be a symplectic manifold. Define a mapping, called the \bld{Poisson bracket on the algebra of classical observables}\index{Poisson!bracket},
	\begin{equation*}
		\cbr{\cdot,\cdot} : C^\infty(M) \times C^\infty(M) \to C^\infty(M)
	\end{equation*}
	\noindent by
	\begin{equation*}
		\cbr{f,g} := \omega(X_f,X_g).
	\end{equation*}
\end{definition}

\begin{definition}[{Poisson Algebra, \cite[Definition~18.7]{silva:sg:2008}}]
	A \bld{Poisson algebra} is defined to be a real commutative algebra $\mathfrak{p}$ together with a Lie bracket $\{\cdot,\cdot\}$ on $\mathfrak{p}$ satisfying the \bld{Leibniz rule}
	\begin{equation*}
		\{f,gh\} = h\{f,g\} + g\{f,h\} \qquad \forall f,g,h \in \mathfrak{p}.
	\end{equation*} 
\end{definition}

\begin{lemma}
	Let $(M,\omega)$ be a symplectic manifold. Then $(C^\infty(M),\{\cdot,\cdot\})$ is a Poisson algebra.
\end{lemma}

\begin{proof}
	The bilinearity and antisymmetry of the Poisson bracket is immediate from the definition. Moreover, the Leibniz rule follows from the computation
	\begin{align*}
		\{f,gh\} &= \omega(X_f,X_{gh})\\
		&= \omega(X_f, h X_g + g X_h)\\
		&= h \omega(X_f,X_g) + g \omega(X_f,X_h)\\
		&= h\{f,g\} + g \{f,h\}
	\end{align*} 
	\noindent for all $f,g,h \in C^\infty(M)$. For proving the Jacobi identity, we claim that
	\begin{equation}
        \label{eq:Jacobi}
		X_{\{f,g\}} = [X_f,X_g] \qquad \forall f,g \in C^\infty(M).
	\end{equation}
    Indeed, we compute
    \begin{align*}
		i_{[X_f,X_g]}\omega &= L_{X_f}i_{X_g}\omega - i_{X_g}L_{X_f}\omega\\
		&= -L_{X_f}dg - i_{X_g}(i_{X_f}d\omega + di_{X_f}\omega)\\
		&= -(i_{X_f}ddg + di_{X_f}dg) + i_{X_g}(ddf)\\
		&= -di_{X_f}dg\\
		&= di_{X_f}i_{X_g}\omega\\
		&= -d\{f,g\}. 
	\end{align*}
    Using formula \eqref{eq:Jacobi}, we compute
	\allowdisplaybreaks
	\begin{align*}
		\{f,\{g,h\}\} &= \omega(X_f,X_{\{g,h\}})\\
		&= \omega(X_f,[X_g,X_h])\\
		&= i_{X_f}\omega[X_g,X_h]\\
		&= -df[X_g,X_h]\\
		&= -[X_g,X_h]f\\
		&= -X_gX_hf + X_hX_gf\\
		&= -X_g\{h,f\} + X_h\{g,f\}\\
		&= -\{g,\{h,f\}\} + \{h,\{g,f\}\}\\
		&= -\{g,\{h,f\}\} - \{h,\{f,g\}\}.
	\end{align*}
	{}
\end{proof}

\begin{corollary}
	Let $(M,\omega)$ be a symplectic manifold. Then
	\begin{equation*}
		(C^\infty(M),\{\cdot,\cdot\}) \to (\mathfrak{X}(M),[\cdot,\cdot]), \qquad f \mapsto X_f
	\end{equation*}
	\noindent is a Lie algebra homomorphism.
\end{corollary}

\begin{proof}
    This immediately follows from formula \eqref{eq:Jacobi}.
\end{proof}
    
\begin{corollary}
    Let $(M,\omega)$ be a symplectic manifold and $\varphi \in \Symp(M,\omega)$. Then 
	\begin{equation*}
		(C^\infty(M),\{\cdot,\cdot\}) \to (C^\infty(M),\{\cdot,\cdot\}), \qquad f \mapsto \varphi^*f
	\end{equation*}
	\noindent is a Lie algebra homomorphism.
\end{corollary}
	
\begin{proof}
    This follows from Lemma \ref{lem:Jacobi}.
\end{proof}

\begin{lemma}[{Evolution Equation, \cite[Corollary~3.3.15]{abrahammarsden:cm:1978}}]
	\label{lem:evolution_equation}
	Let $(M,\omega,H)$ be a Hamiltonian system. Then
	\begin{equation*}
		\frac{d}{dt}f \circ \phi^{X_H}_t = \{H,f\} \circ \phi^{X_H}_t \qquad \forall f \in C^\infty(M),	
	\end{equation*}
	\noindent whenever either side is defined.
\end{lemma}

\begin{proof}
	Using Fisherman's formula \ref{lem:fisherman} we compute
	\begin{align*}
		\frac{d}{dt}f \circ \phi^{X_H}_t &= \frac{d}{dt} \del[1]{\phi_t^{X_H}}^*f\\
		&= \del[1]{\phi_t^{X_H}}^* L_{X_H} f\\
		&= \del[1]{\phi_t^{X_H}}^*\{H,f\}\\
		&= \{H,f\} \circ \phi^{X_H}_t
	\end{align*}
	\noindent for all $f \in C^\infty(M)$.
\end{proof}

\begin{remark}[{The Process of Measurement, \cite[Section~2.8]{takhtajan:qm:2008}}]
	The evolution equation \eqref{lem:evolution_equation} can be stated more concisely as
	\begin{equation*}
		\frac{df}{dt} = \{H,f\} \qquad \forall f \in C^\infty(M),
	\end{equation*}
	\noindent and thus dictates the time-evolution of a classical observable in a Hamiltonian system. If we denote by $\mathscr{M}(M)$ the convex space of all probability measures on $M$, then the above evolution equation together with
	\begin{equation*}
		\frac{d\mu}{dt} = 0 \qquad \forall \mu \in \mathscr{M}(M)
	\end{equation*}
	\noindent constitutes \bld{Hamilton's description of classical mechanics}. Hamilton's picture is commonly used for mechanical systems consisting of few interacting particles. There is another description, called \bld{Liouville's description of classical mechanics}, where classical observables do not depend on time but there is an analog evolution equation for the probability measures on $M$. This picture is commonly used in classical statistical mechanics, where we have macroscopic systems, that is, classical systems with many interacting particles.
\end{remark}

\begin{corollary}[{Preservation of Energy, \cite[Theorem~2.2.2]{frauenfelderkoert:3bp:2018}}]
	\label{cor:preservation_of_energy}
	Let $(M,\omega,H)$ be a Hamiltonian system. Then
	\begin{equation*}
		H\del[1]{\phi^{X_H}_t(x)} = H(x) \qquad \forall x \in M, 
	\end{equation*}
	\noindent whenever the left side is defined.
\end{corollary}

\begin{proof}
	Using Lemma \ref{lem:evolution_equation} we compute
	\begin{equation*}
		\frac{d}{dt}H \circ \phi^{X_H}_t = \{H,H\} \circ \phi^{X_H}_t = 0
	\end{equation*}
	\noindent by antisymmetry of the Poisson bracket.
\end{proof}

Corollary \ref{cor:preservation_of_energy} has many important applications.

\begin{figure}[h!tb]
	\centering
	\includegraphics[width=.85\textwidth]{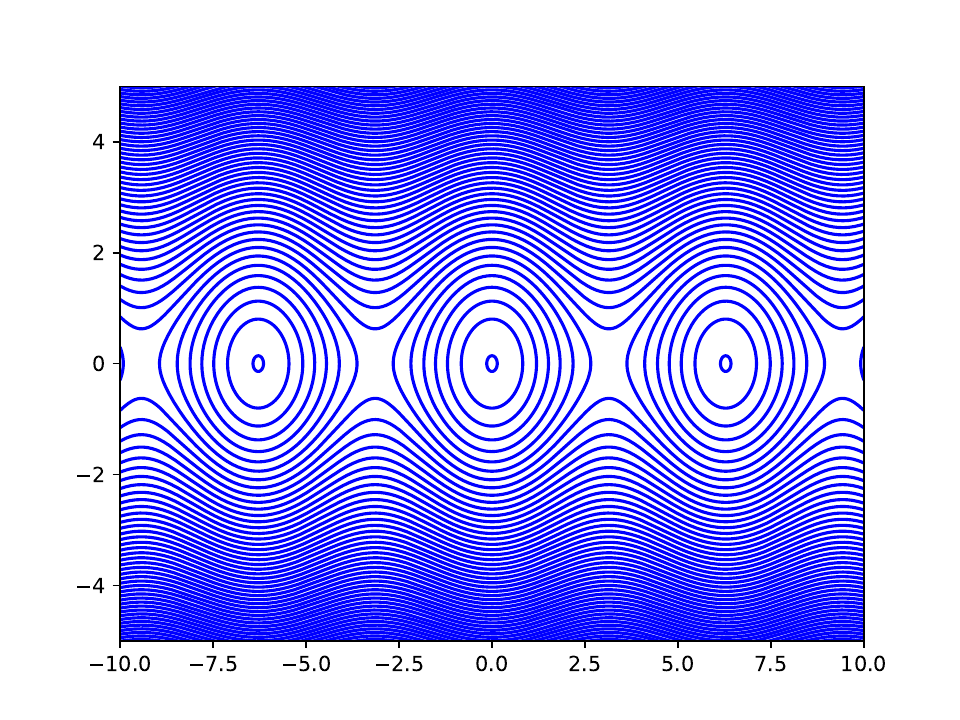}
	\caption{Some level sets $H^{-1}(c)$, $c \geq 0$, of the mathematical pendulum}
	\label{fig:level_sets_pendulum}
\end{figure}

\begin{example}[The Mathematical Pendulum]
    \label{ex:pendulum}
	Consider the mechanical Hamiltonian system $(T^*\mathbb{R},dp \wedge dq, H)$, where
\begin{equation*}
	H \colon \mathbb{R} \times \mathbb{R} \to \mathbb{R}, \qquad H(q,p) := \frac{1}{2}p^2 + (1 - \cos(q)) 
\end{equation*}
\noindent models the motion of a frictionless pendulum. By Preservation of Energy \ref{cor:preservation_of_energy}, we conclude that the integral curves of $X_H$ are contained in the level sets $H^{-1}(c)$, $c \geq 0$. Some of these level sets are depicted in Figure \ref{fig:level_sets_pendulum}. Is there a more formal way to describe these level curves? We follow \cite[Introduction]{cushman:cm:2015}. One can check that
\begin{equation*}
	\Crit H = \{(k\pi,0) : k \in \mathbb{Z}\},
\end{equation*}
\noindent and that $\Hess H(q,p)$ is invertible for every $(q,p) \in \Crit H$. In particular, $\ker \Hess H(q,p)$ is the trivial vector space and thus $H$ is a Morse function. If $k$ is even, then the Morse index of $(k\pi,0)$ is zero and thus the critical points $(k\pi,0)$ are relative minima. If $k$ is odd, then the Morse index of $(k\pi,0)$ is one, and thus the critical points $(k\pi,0)$ are saddle points. By the Morse lemma \cite[Lemma~3.11]{banyagahurtubise:mh:2004}, we can write
\begin{equation*} 
    H(q,p) = f(2k\pi,0) + \frac{1}{2}\Hess H_{(2k\pi,0)}\del[1]{(q,p),(q,p)} = \frac{1}{2}\del[1]{q^2 + p^2}.
\end{equation*}

\begin{exercise}
    How does $H$ look like near the critical points $\del[1]{(2k + 1)\pi,0}$ for $k \in \mathbb{Z}$?
\end{exercise}

\begin{solution}
    By the Morse lemma we get
    \begin{equation*} 
        H(q,p) = f\del[1]{(2k + 1)\pi,0} + \frac{1}{2}\Hess H_{((2k + 1)\pi,0)}\del[1]{(q,p),(q,p)} = \frac{1}{2}\del[1]{p^2 - q^2}.
\end{equation*}
Thus we get hyperbolic paraboloids.
\end{solution}

For a more extensive treatment of the mathematical pendulum in the spirit of dynamical systems, see the excellent book \cite[Section~III.4]{zehnder:ds:2010}.
\end{example}

\begin{exercise}
    Compute the equation of motion for the mathematical pendulum using
    \begin{enumerate}[label=\alph*.]
        \item Hamilton's equations.
        \item Newton's second law.
    \end{enumerate}
    Compare these two methods. Which method is easier to handle?
\end{exercise}

\begin{solution}
    \mbox{}
    \begin{enumerate}[label=\alph*.]
        \item Using formula \eqref{eq:hamiltonian_equations}, we get
        \begin{equation*}
            \dot{q} = p \qquad \text{and} \qquad \dot{p} = -\sin(q).
        \end{equation*}
        Combining these equations yields
        \begin{equation*}
            \ddot{q} = -\sin(q).
        \end{equation*}
        \item The resulting force is given as the tangential component of the gravitational force, and thus $F = -\sin(q)$. Using Newton's second law yields
        \begin{equation*}
            \ddot{q} = F = - \sin(q).
        \end{equation*}
    \end{enumerate}
    In this case, both methods are relatively easy. However, computing the resulting force in a general physical system is usually hard, and in some cases as the brachistochrone problem practically impossible. The drawback of the Hamiltonian equations is, that there is no way to model for example a mathematical pendulum with friction and deriving the equations from a variational principle is conceptually challenging. Indeed, the equation of motion modelling linear friction is given by
        \begin{equation}
            \label{eq:perturbed}
            \ddot{q} = -\mu \dot{q} -\sin q, \qquad \mu > 0.
        \end{equation}
        This equation is equivalent to the first order system
        \begin{equation*}
            \begin{pmatrix}
                \dot{q}\\
                \dot{p}
            \end{pmatrix} = \begin{pmatrix}
                p\\
                -\mu p - \sin(q)
            \end{pmatrix}.
        \end{equation*}
        This system yields the Hamiltonian function
        \begin{equation*}
            H \colon \mathbb{R}^2 \to \mathbb{R}, \qquad H(q,p) = \frac{1}{2}p^2 + \mu qp - \cos(q),
        \end{equation*}
        \noindent which corresponding Hamiltonian equations do not give rise to \eqref{eq:perturbed}. Even introducing a time-dependent Hamiltonian function
        \begin{equation*}
            H \colon \mathbb{R}^2 \times \mathbb{R} \to \mathbb{R}, \qquad H(q,p) := e^{-\mu t}\del[3]{\frac{1}{2}p^2 + (1 - \cos (q))}
        \end{equation*}
        \noindent does not give the right equation of motion as
        \begin{equation*}
            \dot{q} = \frac{\partial H}{\partial p} = e^{-\mu t}p \qquad \text{and} \qquad \dot{p} = -\frac{\partial H}{\partial q} = -e^{-\mu t}\sin(q)
        \end{equation*}
        \noindent implies
        \begin{equation*}
            \ddot{q} = -\mu e^{-\mu t} p + e^{-\mu t} \dot{p} = -\mu \dot{q} - e^{-2\mu t} \sin(q).
        \end{equation*}
\end{solution}

\subsection{Noether's Theorem}
Motivated by Lemma \ref{lem:evolution_equation}, we give the following definition of a preserved quantity in a Hamiltonian system.

\begin{definition}[Integral of Motion]
	An \bld{integral of motion} for a Hamiltonian system $(M,\omega,H)$ is defined to be a smooth function $I \in C^\infty(M)$ such that $\{H,I\}= 0$.
\end{definition}

\begin{remark}
	It is easy to check that the integrals of motion form a Lie subalgebra of the Poisson algebra of classical observables $(C^\infty(M),\{\cdot,\cdot\})$.
\end{remark}

Let $\psi \colon G \times M \to M$ be a smooth left action of a Lie group $G$ on a smooth manifold $M$ and denote by $\mathfrak{g} := T_e G \cong \mathfrak{X}_L(G)$ the corresponding Lie algebra. Every $\xi \in \mathfrak{g}$ determines a smooth global flow on $M$ by $(t,x) \mapsto \psi_{\exp(-t\xi)}(x)$, where
\begin{equation*}
	\exp \colon \mathfrak{g} \to G, \qquad \exp(\xi) := \gamma_\xi(1)
\end{equation*}
\noindent denotes the exponential map and $\gamma_\xi$ is the integral curve of the left-invariant vector field $X_\xi$ starting at $e$ with $\dot{\gamma}_\xi(0) = \xi$. Note that if there exists a bi-invariant Riemannian metric on $G$, then this exponential map coincides with the exponential map of the associated Levi--Civita connection at the identity \cite[Problem~5--8.~(c)]{lee:dg:2018}. Define $\widehat{\xi} \in \mathfrak{X}(M)$ to be the infinitesimal generator of this flow, that is,
\begin{equation*}
	\widehat{\xi}_x = \frac{d}{dt}\bigg\vert_{t = 0} \psi_{\exp(-t\xi)}(x) \qquad \forall x \in M.
\end{equation*}
Then
\begin{equation*}
	(\mathfrak{g},[\cdot,\cdot]) \to (\mathfrak{X}(M),[\cdot,\cdot]), \qquad \xi \mapsto \widehat{\xi}
\end{equation*}
\noindent is a Lie algebra homomorphism.

\begin{definition}[Weakly Hamiltonian Action]
	Let $(M,\omega)$ be a symplectic manifold. A smooth action $\psi \colon G \times M \to M$ of a Lie group $G$ is said to be \bld{weakly Hamiltonian}, if the map $\psi_g \in \Diff(M)$ is a symplectomorphism for all $g \in G$ and there exists a linear map
	\begin{equation*}
		\mu \colon \mathfrak{g} \to C^\infty(M),
	\end{equation*}
	\noindent called a \bld{momentum map}, such that the diagram
	\begin{equation*}
		\begin{tikzcd}
			C^\infty(M) \arrow[rr,"f \mapsto X_f"] & & \mathfrak{X}(M,\omega)\\
			& \mathfrak{g} \arrow[ul,"\mu"]\arrow[ur,"\xi \mapsto \widehat{\xi}"']
		\end{tikzcd}
	\end{equation*}
	\noindent commutes.
\end{definition}

\begin{definition}
	A weakly Hamiltonian action of a Lie group $G$ on a symplectic manifold $(M,\omega)$ is called 
	\begin{itemize}
		\item \bld{Hamiltonian}, if the momentum map $\mu \colon \mathfrak{g} \to C^\infty(M)$ is $G$-equivariant with respect to the adjoint action of $G$ on its associated Lie algebra $\mathfrak{g}$ and the induced action of $G$ on $C^\infty(M)$, that is
			\begin{equation*}
				\mu(\Ad_{g^{-1}}(\xi)) = \mu(\xi) \circ \psi_g \qquad \forall g \in G, \xi \in \mathfrak{g},
			\end{equation*}
			\noindent where 
			\begin{equation*}
				\Ad_{g^{-1}}(\xi) := \frac{d}{dt}\bigg\vert_{t = 0} g^{-1}\exp(t\xi)g.
			\end{equation*}
		\item \bld{Poisson}, if the associated momentum map
			\begin{equation*}
				\mu \colon (\mathfrak{g},[\cdot,\cdot]) \to (C^\infty(M),\{\cdot,\cdot\})
			\end{equation*}
			\noindent is a Lie algebra homomorphism.
	\end{itemize}
\end{definition}

\begin{remark}
    By \cite[p.~33--34]{frauenfelderkoert:3bp:2018}, for a connected Lie group, a weakly Hamiltonian action is Hamiltonian if and only if it is Poisson.
\end{remark}

For showing existence and uniqueness results for Poisson actions, we recall the basic notions of Lie algebra cohomology.	Let $\mathfrak{g}$ be a Lie algebra. Define
	\begin{equation*}
		C^k := \Lambda^k \mathfrak{g}^*
	\end{equation*}
	\noindent and $d \colon C^k \to C^{k + 1}$ by
	\begin{equation*}
		df(X_0,\dots,X_k) := \sum_{0 \leq i < j \leq k} (-1)^{i + j} f([X_i,X_j],X_0,\dots,\overline{X}_i,\dots,\overline{X}_j,\dots,X_k).
	\end{equation*}
	Then one checks that $d \circ d = 0$. The resulting nonnegative chain complex is called the \bld{Chevalley--Eilenberg cochain complex}. Then the \bld{$k$-th cohomology group of $\mathfrak{g}$} is defined by
		\begin{equation*}
			H^k(\mathfrak{g};\mathbb{R}) := \frac{\ker d \colon C^k \to C^{k + 1}}{\im d \colon C^{k - 1} \to C^k}.
		\end{equation*}
  By \cite[Theorem~26.1]{cannasdasilva:sg:2008}, we have that $H^*(\mathfrak{g};\mathbb{R}) \cong H^*_{\operatorname{dR}}(G)$ for the Lie algebra $\mathfrak{g}$ of a compact and connected Lie group $G$.

  \begin{exercise}[$H^3(\mathfrak{sl}_2(\mathbb{R}),\mathbb{R})$]
      Consider the \bld{special linear Lie algebra} $\mathfrak{sl}_2(\mathbb{R},\mathbb{R})$ with ordered basis 
      	\begin{equation*}
		e_1 := \begin{pmatrix}
			1 & 0\\
			0 & -1
		\end{pmatrix} \qquad
		e_2 := \begin{pmatrix}
			0 & 1\\
			0 & 0
		\end{pmatrix} \qquad
		e_3 := \begin{pmatrix}
			0 & 0\\
			1 & 0
		\end{pmatrix}.
	\end{equation*}
    Show that $H^3(\mathfrak{sl}_2(\mathbb{R}),\mathbb{R}) \cong \mathbb{R}$.
  \end{exercise}

  \begin{solution}
      The crucial portion of the Eilenberg--Chevalley cochain complex is given by
	\begin{equation*}
		\begin{tikzcd}
			\dots \arrow[r] & \Lambda^2\mathfrak{sl}_2(\mathbb{R}) \arrow[r,"d"] & \Lambda^3\mathfrak{sl}_2(\mathbb{R}) \arrow[r] & 0.
		\end{tikzcd}
	\end{equation*}
	We compute
	\begin{align*}
		df(e_1,e_2,e_3) =& e_1f(e_2,e_3) - e_2f(e_1,e_3) + e_3f(e_1,e_2)\\
		&- f([e_1,e_2],e_3) + f([e_1,e_3],e_2) - f([e_2,e_3],e_1)\\
		=& -2f(e_2,e_3) - 2f(e_3,e_2) - f(e_1,e_1)\\
		=& -2f(e_2,e_3) + 2f(e_2,e_3)\\
		=& 0.
	\end{align*}
	 It follows that
	\begin{equation*}
		H^3(\mathfrak{sl}_2(\mathbb{R});\mathbb{R}) \cong \Lambda^3 \mathfrak{sl}_2(\mathbb{R}) \cong \mathbb{R}.
	\end{equation*}
  \end{solution}
  
\begin{theorem}[Uniqueness of Momentum Maps for Poisson Actions]
	\label{thm:uniqueness_momentum_map}
	Let $\mu$ and $\tilde{\mu}$ be two momentum maps for a Poisson $G$-action on a connected symplectic manifold. If for the first Lie algebra cohomology group $H^1(\mathfrak{g};\mathbb{R}) = 0$ holds, then $\mu = \tilde{\mu}$.
\end{theorem}

\begin{proof}
	By assumption there exists $\sigma \in \mathfrak{g}^*$ such that
	\begin{equation*}
		\mu(\xi) - \tilde{\mu}(\xi) = \sigma(\xi) \qquad \forall \xi \in \mathfrak{g}.
	\end{equation*}
	Since both $\mu$ and $\tilde{\mu}$ are Lie algebra homomorphisms, we have that $d\sigma = 0$. Indeed, for $\xi,\eta \in \mathfrak{g}$ we compute
	\begin{equation*}
		d\sigma(\xi,\eta) = \sigma([\eta,\xi]) = \mu([\eta,\xi]) - \tilde{\mu}([\eta,\xi]) = \{\mu(\eta),\mu(\xi)\} - \{\tilde{\mu}(\eta),\tilde{\mu}(\xi)\} = 0.
	\end{equation*}
	Thus $\sigma \in H^1(\mathfrak{g};\mathbb{R}) = 0$, implying $\sigma = 0$ and the statement follows.
\end{proof}

\begin{theorem}[Existence of Poisson Actions]
	\label{thm:existence_Poisson_action}
	Suppose we are given a weakly Hamiltonian $G$-action on a connected symplectic manifold $(M,\omega)$. If $H^2(\mathfrak{g};\mathbb{R}) = 0$, then the action is Poisson.
\end{theorem}

\begin{proof}
	For $\xi,\eta \in \mathfrak{g}$ we compute
	\begin{equation*}
		X_{\mu([\xi,\eta])} = \widehat{[\xi,\eta]} = [\widehat{\xi},\widehat{\eta}] = [X_{\mu(\xi)},X_{\mu(\eta)}] = X_{\{\mu(\xi),\mu(\eta)\}}
	\end{equation*}
	\noindent using Proposition \eqref{eq:Jacobi}. Thus by connectedness of $M$ there exists $\tau \in \Lambda^2\mathfrak{g}^*$ such that
	\begin{equation*}
		\{\mu(\xi),\mu(\eta)\} - \mu([\xi,\eta]) = \tau(\xi,\eta) \qquad \forall \xi,\eta \in \mathfrak{g}.
	\end{equation*}
	Invoking the Jacobi identity for the Lie as well as the Poisson bracket, yields $d\tau = 0$ and so $\tau \in H^2(\mathfrak{g};\mathbb{R}) = 0$. Hence there exists $\sigma \in H^1(\mathfrak{g};\mathbb{R})$ such that $\tau = d\sigma$. The momentum map
	\begin{equation*}
		\mathfrak{g} \to C^\infty(M), \qquad \xi \mapsto \mu(\xi) - \sigma(\xi)
	\end{equation*}
	\noindent is a Lie algebra homomorphism. 
\end{proof}

Recall, that a Lie algebra $\mathfrak{g}$ is said to be semisimple if $\mathfrak{g}$ does not admit any nontrivial abelian ideals. A Lie group $G$ is called semisimple, if its associated Lie algebra $\mathfrak{g} := T_e G$ is semisimple.

\begin{corollary}
	Let $G$ be a semisimple Lie group. Then every weakly Hamiltonian $G$-action on a connected symplectic manifold is Poisson and admits a unique momentum map.
\end{corollary}

\begin{proof}
	The statement immediately follows from Theorem \ref{thm:uniqueness_momentum_map} and \ref{thm:existence_Poisson_action} as 
	\begin{equation*}
		H^1(\mathfrak{g};\mathbb{R}) = H^2(\mathfrak{g};\mathbb{R}) = 0
	\end{equation*}
	\noindent by \cite[Theorem~21.1]{chevalleyeilenberg:lie:1948}.
\end{proof}

\begin{remark}
    From \cite[Theorem~21.1]{chevalleyeilenberg:lie:1948} it immediately follows that $\mathbb{S}^n$ is a Lie group for $n = 1$ and $n = 3$ only \cite[Theorem~21.3]{chevalleyeilenberg:lie:1948}.
\end{remark}

\begin{lemma}[{Momentum Lemma, \cite[Exercise~5.2.2]{mcduffsalamon:st:2017}}]
	\label{lem:momentum_lemma}
	Let $\psi \colon G \times M \to M$ be a smooth Lie group action on an exact symplectic manifold $(M,\lambda)$ such that $\psi^*_g \lambda = \lambda$ for all $g \in G$ holds. Then the action is Hamiltonian and Poisson with momentum map
	\begin{equation*}
		\mu(\xi) = i_{\widehat{\xi}}\lambda, \qquad \forall \xi \in \mathfrak{g}.
	\end{equation*}
\end{lemma}

\begin{proof}
	We show the result in four steps. Obviously, $\psi_g^*d\lambda = d\lambda$ for all $g \in G$.
	
	\emph{Step 1: $\psi$ is a weakly Hamiltonian action.} Let $\xi \in \mathfrak{g}$. We compute 
			\allowdisplaybreaks
			\begin{align*}
				i_{\widehat{\xi}}d\lambda &= L_{\widehat{\xi}}\lambda - di_{\widehat{\xi}}\lambda\\
				&= \frac{d}{dt}\bigg\vert_{t = 0} \psi^*_{\exp(-t\xi)}\lambda - d\mu(\xi)\\
				&= \frac{d}{dt}\bigg\vert_{t = 0} \lambda - d\mu(\xi)\\
				&= -d\mu(\xi).
			\end{align*}
		
			\emph{Step 2: $\psi^*_g\widehat{\xi} = \widehat{\Ad_{g^{-1}}(\xi)}$ for all $g \in G$ and $\xi \in \mathfrak{g}$.} We have a commutative diagram
			\begin{equation*}
				\begin{tikzcd}
					\mathfrak{g} \arrow[rr,"\Ad_{g^{-1}}"]\arrow[dd,"\exp"'] & & \mathfrak{g}\arrow[dd,"\exp"]\\
					\\
					G \arrow[rr,"C_{g^{-1}}"'] & &  G,\\
				\end{tikzcd}
			\end{equation*}
			\noindent where $C_{g^{-1}}(h) = g^{-1}hg$ denotes the conjugation action on $G$. Let $g \in G$. Then we compute
			\allowdisplaybreaks
			\begin{align*}
				\psi^*_g \widehat{\xi} &= D\psi_{g^{-1}} \circ \widehat{\xi} \circ \psi_g\\
				&= D\psi_{g^{-1}} \circ \frac{d}{dt}\bigg\vert_{t = 0} \psi_{\exp(-t\xi)} \circ \psi_g\\
				&= \frac{d}{dt}\bigg\vert_{t = 0} \psi_{g^{-1}} \circ \psi_{\exp(-t\xi)} \circ \psi_g\\
				&= \frac{d}{dt}\bigg\vert_{t = 0} \psi_{g^{-1}\exp(-t\xi)g}\\
				&= \frac{d}{dt}\bigg\vert_{t = 0} \psi_{\exp(-t\Ad_{g^{-1}}(\xi))}\\
				&= \widehat{\Ad_{g^{-1}}(\xi)}.
			\end{align*}
		
			\emph{Step 3: $\psi$ is a Hamiltonian action.} Using Step 2, we compute
			\begin{align*}
				\mu(\Ad_{g^{-1}}(\xi)) &= i_{\widehat{\Ad_{g^{-1}}(\xi)}}\lambda\\
				&= i_{\psi_g^* \widehat{\xi}}\lambda\\
				&= \lambda(D\psi_{g^{-1}} \circ \widehat{\xi} \circ \psi_g)\\
				&= \psi_{-g}^*\lambda(\widehat{\xi} \circ \psi_g)\\
				&= i_{\widehat{\xi}}\lambda \circ \psi_g\\
				&= \mu(\xi) \circ \psi_g
			\end{align*}
			\noindent for all $g \in G$ and $\xi \in \mathfrak{g}$.
		
			\emph{Step 4: $\psi$ is a Poisson action.} For $\xi,\eta \in \mathfrak{g}$ we compute
			\begin{align*}
				\mu[\xi,\eta] &= i_{\widehat{[\xi,\eta]}}\lambda\\
				&= i_{[\widehat{\xi},\widehat{\eta}]}\lambda\\
				&= L_{\widehat{\xi}}i_{\widehat{\eta}}\lambda - i_{\widehat{\eta}}L_{\widehat{\xi}}\lambda\\
				&= L_{\widehat{\xi}}i_{\widehat{\eta}}\lambda\\
				&= L_{\widehat{\xi}}\mu(\eta)\\
				&= \widehat{\xi}\mu(\eta)\\
				&= X_{\mu(\xi)}\mu(\eta)\\
				&= \{\mu(\xi),\mu(\eta)\}.
			\end{align*}
	{}
\end{proof}

\begin{example}[Cotangent Lift]
    \label{ex:cotangent_lift}
	Let $M$ be a smooth manifold and $\varphi \in \Diff(M)$. Define a map $D\varphi^\dagger \colon T^* M \to T^* M$, called the \bld{cotangent lift of $\varphi$}, by
	\begin{equation*}
		D\varphi^\dagger(q,p)(v) := p\del[1]{D\varphi^{-1}(v)} \qquad \forall v \in T_{\varphi(q)}M.
	\end{equation*}
	Then one checks that $(D\varphi^\dagger)^*\lambda = \lambda$. Thus if we have a smooth Lie group action $\psi \colon G \times M \to M$, we get an induced action
	\begin{equation*}
		\Psi \colon G \times T^* M \to T^* M, \qquad \Psi_g := D\psi_g^\dagger.
	\end{equation*}
	By the Momentum Lemma \ref{lem:momentum_lemma} we conclude that $\Psi$ is Hamiltonian and Poisson with momentum map
	\allowdisplaybreaks
	\begin{align*}
		\mu(\xi)(q,p) &= \lambda_{(q,p)}\del[1]{\widehat{\xi}(q,p)}\\
		&= p\del[1]{D\pi\del[1]{\widehat{\xi}(q,p)}}\\
		&= p\del[3]{\frac{d}{dt}\bigg\vert_{t = 0}\pi\del[1]{\Psi_{\exp(-t\xi)}(q,p)}}\\
		&= p\del[3]{\frac{d}{dt}\bigg\vert_{t = 0}\psi_{\exp(-t\xi)}(\pi(q,p))}\\
		&= p\del[3]{\frac{d}{dt}\bigg\vert_{t = 0}\psi_{\exp(-t\xi)}(q)}
	\end{align*}
	\noindent for all $(q,p) \in T^* M$.
\end{example}

\begin{definition}[Symmetry Group]
	A Lie group $G$ is said to be a \bld{symmetry group of a Hamiltonian system $(M,\omega,H)$}\index{Symmetry group}, if there exists a weakly Hamiltonian action $\psi$ of $G$ on $(M,\omega)$, such that $H \circ \psi_g = H$ for all $g \in G$.
\end{definition}

\begin{lemma}[Noether's Theorem]
	\label{lem:Noethers_theorem}
	Let $G$ be a symmetry group of a Hamiltonian system $(M,\omega,H)$. Then $\mu(\xi)$ is an integral of motion for all $\xi \in \mathfrak{g}$.
\end{lemma}

\begin{proof}
	For $\xi \in \mathfrak{g}$ we compute
	\begin{equation*}
		\{\mu(\xi),H\} = X_{\mu(\xi)}H = \widehat{\xi}H = \frac{d}{dt}\bigg\vert_{t = 0}H \circ \psi_{\exp(-t\xi)} = \frac{d}{dt}\bigg\vert_{t = 0}H = 0.
	\end{equation*}
	{}
\end{proof}

\begin{example}[The Kepler Problem]
    \label{ex:kepler}
    The Hamiltonian function of the $n$-dimensional Kepler problem is given by
    \begin{equation}
        \label{eq:kepler}
        H \colon T^*(\mathbb{R}^n \setminus \{0\}) \to \mathbb{R}, \qquad H(q,p) := \frac{1}{2}\norm[0]{p}^2 - \frac{1}{\norm[0]{q}}.
    \end{equation}
    Define an $\operatorname{SO}(n)$-action on $\mathbb{R}^* \setminus \{0\}$ by
    \begin{equation*}
        \psi\colon \operatorname{SO}(n) \times \mathbb{R}^* \setminus \{0\} \to \mathbb{R}^* \setminus \{0\}, \qquad \psi_A(q) := Aq.
    \end{equation*}
    Then by Example \ref{ex:cotangent_lift}, we get an induced $\operatorname{SO}(n)$-action on 
    \begin{equation*}
        T^*(\mathbb{R}^n \setminus \{0\}) \cong (\mathbb{R}^n \setminus \{0\}) \times \mathbb{R}^n
    \end{equation*}
    \noindent with momentum map
    \begin{equation*}
        \mu \colon \mathfrak{so}(n) \to C^\infty((\mathbb{R}^n \setminus \{0\}) \times \mathbb{R}^n), \qquad \mu(B)(q,p) = q^TBp,
    \end{equation*}
    \noindent where 
    \begin{equation*}
        \mathfrak{so}(n) = \cbr[1]{A \in \mathbb{R}^{n \times n} : A^T = -A}.
    \end{equation*}
\end{example}

\begin{exercise}
    Write down the three integrals of motion in the Kepler problem given by Noether's theorem \ref{lem:Noethers_theorem} for the case $n = 3$. What is the physical interpretation of these integrals of motion?
\end{exercise}

\begin{solution}
    In the spatial case $n = 3$ we have that
    \begin{equation*}
        e_1 := \begin{pmatrix}
            0 & 0 & 0\\
            0 & 0 & 1\\
            0 & -1 & 0
        \end{pmatrix}, \qquad e_2 := \begin{pmatrix}
            0 & 0 & 1\\
            0 & 0 & 0\\
            -1 & 0 & 0
        \end{pmatrix} \qquad \text{and} \qquad e_3 := \begin{pmatrix}
            0 & 1 & 0\\
            -1 & 0 & 0\\
            0 & 0 & 0
        \end{pmatrix}
    \end{equation*}
    \noindent is a basis for $\mathfrak{so}(3)$. The integrals of motion $\mu(e_j)$ are precisely the components $L_j$ of the \bld{angular momentum}
    \begin{equation*}
        L \colon T^*(\mathbb{R}^n \setminus \{0\}) \to \mathbb{R}^3, \qquad L(q,p) := q \times p. 
    \end{equation*}
\end{solution}

For a much more detailed analysis of the Kepler problem as well its connection to the restricted three body problem see \cite{frauenfelderkoert:3bp:2018} and the excellent survey article \cite{moreno:3bp:2022}

\begin{exercise}[Summary Hamiltonian Formalism]
    \mbox{}
    \begin{enumerate}[label=\alph*.]
    \item Write a summary of the important concepts of this section about the Hamiltonian formalism and \emph{think} of one answer to each of the following questions:
        \begin{itemize}
            \item Which concepts were new to you?
            \item Which concepts you were already familiar with?
            \item Which concepts were hard to understand and should be elaborated more on?
        \end{itemize}
        \item \emph{Pair} up with another student and discuss your answers.
        \item \emph{Share} your results with the other students and the lecturer.
    \end{enumerate}
\end{exercise}

\begin{summary*}{The Hamiltonian Formalism}
    \begin{itemize}
        \item A Hamiltonian system is a triple $(M,\omega,H)$, where $(M,\omega)$ is a symplectic manifold, and $H \in C^\infty(M)$ is a classical observable.
        \item The dynamics of a Hamiltonian system is governed by the evolution equation
        \begin{equation*}
            \frac{d}{dt} f \circ \phi_t = \{H,f\} \circ \phi_t,
        \end{equation*}
        \noindent where $\phi$ denotes the flow of the Hamiltonian vector field $X_H$ and $\{\cdot,\cdot\}$ is the Poisson bracket.
        \item In Darboux coordinates, the integral curves of the Hamiltonian vector field satisfy the Hamiltonian equations
        \begin{equation*}
            \dot{q} = \frac{\partial H}{\partial p} \qquad \text{and} \qquad \dot{p} = -\frac{\partial H}{\partial q}.
        \end{equation*}
        \item If $G$ is a symmetry group, then the integrals of motions are given by $\mu(\xi)$ for every $\xi \in \mathfrak{g} = T_e G$ by Noether's theorem, where 
        \begin{equation*}
            \mu \colon \mathfrak{g} \to C^\infty(M)
        \end{equation*}
        \noindent is the momentum map.
        \item The Hamiltonian formalism is unable to model dissipative systems as the mathematical pendulum with friction. However, the Hamiltonian formalism can be used to analyse problems in celestial mechanics qualitatively, as the Kepler problem or the restricted three body problem. Moreover, the abstract formulation of classical mechanics using the Hamiltonian formalism can also be used in modern physics in the formulation of quantum mechanics.
    \end{itemize}
\end{summary*}
\newpage
\section{The Lagrangian Formalism}
The Lagrangian formalism of classical mechanics is the dual viewpoint of the Hamiltonian formalism. This description initially preceeded the Hamiltonian one and was inspired by variational observations coming from optics (Fermat's principle) and from the solution of the \emph{brachistochrone problem} in 1697. In contrast to Newton's formulation, the Lagrangian formulation is independent of the choice of coordinates and holonomic constraints can easily be implemented.

\subsection{The Legendre Transform}

\begin{definition}[Lagrangian System]
	A \bld{Lagrangian system}\index{Lagrangian!system} is a tuple $(M,L)$, where $M$ is a finite-dimensional smooth manifold, called the \bld{configuration space}, and $L \in C^\infty(TM)$ is a smooth function, called the \bld{Lagrangian function}\index{Lagrangian!function}. Moreover, the tangent bundle $TM$ of the configuration space $M$ is called the \bld{state space}\index{Space!state}.
\end{definition}

\begin{exercise}
    What is the physical interpretation of the standard coordinates $(q,\dot{q})$ on the state space $TM$ of a configuration space $M$?
\end{exercise}

\begin{solution}
    The coordinate $q$ is interpreted as the position and $\dot{q}$ as velocity.
\end{solution}

\begin{remark}
	Infinite-dimensional configuration spaces are treated in classical field theory \cite{diez:field_theory:2014}.
\end{remark}

\begin{remark}
	\label{rem:double_tangent_bundle}
	Let $\pi \colon E \to M$ be a vector bundle, $\mathcal{H}$ a connection on $E$ and $\kappa \colon TE \to E$ be the associated connection map. Then 
	\begin{equation*}
		(D\pi,\kappa) \colon TE \to TM \oplus E
	\end{equation*}
	\noindent is a vector bundle isomorphism along $\pi$. In particular
	\begin{equation*}
		TTM \cong TM \oplus TM
	\end{equation*}
	\noindent as vector bundles for every smooth manifold $M$. Note that this isomorphism depends on the choice of a connection and is therefore not canonical.
\end{remark}

\begin{definition}[Mechanical Lagrangian Function]
    \label{def:mechanical_lagrangian}
	Let $(M,m)$ be a pseudo-Riemannian manifold and $V \in C^\infty(M)$. A \bld{mechanical Lagrangian function} is defined to be the function $L_V \in C^\infty(TM)$
	\begin{equation*}
		L_V(q,\dot{q}) := \frac{1}{2}\norm[0]{\dot{q}}_m^2 - V(q).
	\end{equation*}
\end{definition}

If $F \in C^\infty(M,N)$, then the derivative of $F$ can be interpreted as a vector bundle homomorphism $DF : TM \to F^*TN$. Indeed, define 
\begin{equation*}
	DF(q,\dot{q}) := (q, (F(q),DF_q(\dot{q})))
\end{equation*}
\noindent for any $(q,\dot{q}) \in TM$. If $\pi \colon E \to M$ is a fibre bundle, we can set
\begin{equation*}
	VE := \coprod_{p \in E} \ker D\pi_p.
\end{equation*}
Then $VE$ with the usual footpoint projection is a vector bundle over $E$, called the \bld{vertical bundle of $E$}. Moreover, one can show that $VE$ is isomorphic to $\pi^*E$. Explicitly, the isomorphism $\Phi\colon \pi^*E \to VE$ is given by
\begin{equation}
	\label{eq:pullback_bundle_isomorphism}
	\Phi(\dot{q},v) := \frac{d}{d\varepsilon}\bigg\vert_{\varepsilon = 0}(\dot{q} + \varepsilon v).
\end{equation}
 
\begin{definition}[Noether Form]
	\label{def:Noether_form}
	Let $(M,L)$ be a Lagrangian system. Define the \bld{Noether form}, written $\lambda_L \in \Omega^1(TM)$, by
	\begin{equation}
		\lambda_L(\xi) := dL((\Phi \circ D\pi_{TM})\xi) \qquad \forall \xi \in TTM.
	\end{equation}
\end{definition}

\begin{definition}[Legendre Transform]
	\label{def:Legendre_transform}
	A \bld{Legendre transform of a Lagrangian system $(M,L)$}\index{Transform!Legendre} is defined to be a mapping $\tau_L \in C^\infty(TM,T^*M)$ such that
	\begin{equation*}
		\pi_{T^*M} \circ \tau_L = \pi_{TM} \qquad \text{and} \qquad \tau_L^*\lambda = \lambda_L,
	\end{equation*}
    \noindent where $\lambda \in \Omega^1(T^*M)$ denotes the Liouville form \eqref{eq:liouville}.
\end{definition}

\begin{lemma}[{Fibre derivative, \cite[Definition~3.5.1]{abrahammarsden:cm:1978}}]
    \label{lem:fibre_derivative}
	Let $(M,L)$ be a Lagrangian system. Then $\tau_L \in C^\infty(TM,T^*M)$ is a Legendre transform if and only if
	\begin{equation*}
		\tau_L(\dot{q})(v) = \frac{d}{d\varepsilon}\bigg\vert_{\varepsilon = 0}L_q(\dot{q} + \varepsilon v) \qquad \forall \dot{q}, v \in T_qM, q \in M.
	\end{equation*}
\end{lemma}

\begin{proof}
	Suppose $\tau_L \in C^\infty(TM,T^*M)$ is a Legendre transform. For $\bar{v} \in D\pi_{TM}^{-1}(v)$ we compute on one hand
	\begin{align*}
		(\tau_L^*\lambda)_{(q,\dot{q})}(\bar{v}) &= \lambda_{(q,\tau_L(\dot{q}))}(D\tau_L(\bar{v}))\\
		&= \tau_L(\dot{q})\del[1]{(D\pi_{T^*M} \circ D\tau_L)(\bar{v})}\\
		&= \tau_L(\dot{q})\del[1]{D(\pi_{T^*M} \circ \tau_L)(\bar{v})}\\
		&= \tau_L(\dot{q})\del[1]{D\pi_{TM}(\bar{v})}\\
		&= \tau_L(\dot{q})(v),
	\end{align*}
	\noindent and on the other
	\begin{align*}
		\lambda_L\vert_{(q,\dot{q})}(\bar{v}) &= dL((\Phi \circ D\pi_{TM})(\bar{v}))\\
		&= dL\del[3]{\frac{d}{d\varepsilon}\bigg\vert_{\varepsilon = 0}\del[1]{\dot{q} + \varepsilon D\pi_{TM}(\bar{v})}}\\
		&= \frac{d}{d\varepsilon}\bigg\vert_{\varepsilon = 0}L_q(\dot{q} + \varepsilon D\pi_{TM}(\bar{v}))\\
		&= \frac{d}{d\varepsilon}\bigg\vert_{\varepsilon = 0}L_q(\dot{q} + \varepsilon v).
	\end{align*}
	{}
\end{proof}

\begin{exercise}
    \label{ex:pseudo}
    Compute the Legendre transform for a mechanical Lagrangian function $L_V$ as in Definition \ref{def:mechanical_lagrangian} and show that it is a diffeomorphism. What is the physical interpretation of $L_V$?
\end{exercise}

\begin{solution}
    We make use of the formula given in Lemma \ref{lem:fibre_derivative}. If $m$ is a pseudo-Riemannian metric, we compute
    \begin{equation*}
        \tau_{L_V}(\dot{q})(v) = \frac{1}{2}\frac{d}{d\varepsilon}\bigg\vert_{\varepsilon = 0}\norm[0]{\dot{q} + \varepsilon v}^2 = \frac{1}{2}\frac{d}{d\varepsilon}\bigg\vert_{\varepsilon = 0}m_q(\dot{q} + \varepsilon v,\dot{q} + \varepsilon v) = m_q(\dot{q},v).
    \end{equation*}
    Thus $\tau_L$ is simply the bundle isomorphism
    \begin{equation*}
        TM \to T^*M, \qquad (q,\dot{q}) \mapsto i_{\dot{q}}m.
    \end{equation*}
    The physical interpretation of a mechanical Lagrangian function $L_V$ is kinetic \emph{minus} potential energy.
\end{solution}

\begin{definition}[Energy]
	\label{def:energy}
	The \bld{energy of a Lagrangian system $(M,L)$}\index{Energy!of an autonomous Lagrangian system} is defined to be the function $E_L \in C^\infty(TM)$ given by
	\begin{equation*}
		E_L(q,\dot{q}) := \tau_L(\dot{q})(\dot{q}) - L(q,\dot{q})
	\end{equation*}
	\noindent for $(q,\dot{q}) \in TM$.
\end{definition}

\begin{definition}[Hamiltonian Function]
	\label{def:Hamiltonian_function}
	Let $(M,L)$ be a Lagrangian system such that the Legendre transform $\tau_L$ is a diffeomorphism. The function $H_L \in C^\infty(T^*M)$ defined by
	\begin{equation*}
		H_L := E_L \circ \tau_L^{-1}
	\end{equation*}
	\noindent is called the \bld{Hamiltonian function associated with the Lagrangian function $L$}\index{Hamiltonian!function}.
\end{definition}

\begin{exercise}
    Compute the Hamiltonian function $H_L$ associated with a mechanical Lagrangian function. What is the physical interpretation of $H_L$ and the standard coordinates $(q,p)$ on the phase space $T^*M$? 
\end{exercise}

\begin{solution}
    Using Exercise \ref{ex:pseudo} we compute
    \begin{equation*}
        H_L(q,p) = E_L\del[1]{\tau_L^{-1}(p)} = \frac{1}{2}m_q\del[1]{\tau_L^{-1}(p),\tau_L^{-1}(p)} + V(q) = \frac{1}{2}\norm[0]{p}_{m^*}^2 + V(q).
    \end{equation*}
    The physical interpretation of $H_L$ is kinetic \emph{plus} potential energy, so $H_L$ is just the total energy of the system, with the coordinate $q$ being position and $p$ being momentum.
\end{solution}

\begin{remark}[Tonelli Lagrangians]
	Let $(M^n,L)$ be a Lagrangian system and fix a Riemannian metric $m$ on $M$. The Lagrangian function $L$ is said to be \bld{Tonelli}, if the following conditions are satisfied:
	\begin{enumerate}[label=(T\arabic*),leftmargin=*]
		\item The fibrewise Hessian of $L$ is positive-definite, that is, 
			\begin{equation*}
				\sum_{i,j = 1}^n\frac{\partial^2 L}{\partial \dot{q}_i \partial \dot{q}_j}(q,\dot{q})v_i v_j > 0
			\end{equation*}
			\noindent for all $(q,\dot{q}) \in TM$ and $v := \sum_{i = 1}^n v_i \partial_i \in T_xM$ such that $v \neq 0$.
		\item The Lagrangian function $L$ is fibrewise supercoersive, that is,
			\begin{equation*}
				\lim_{\norm{\dot{q}}_m \to \infty} \frac{L(q,\dot{q})}{\norm{\dot{q}}_m} = +\infty
			\end{equation*}
			\noindent for all $q \in M$.
	\end{enumerate}
	By \cite[Proposition~1.2.1]{mazzucchelli:cm:2012}, for a fibrewise convex Lagrangian function $L$, the associated Legendre transformation $\tau_L \colon TM \to T^*M$ is a diffeomorphism, if and only if $L$ is Tonelli. 
\end{remark}

\begin{exercise}
    Show that every mechanical Lagrangian $L_V$ is Tonelli.
\end{exercise}

\begin{solution}
    \mbox{}
    	\begin{enumerate}[label=(T\arabic*),leftmargin=*]
		\item We compute 
			\begin{equation*}
				\sum_{i,j = 1}^n\frac{\partial^2 L}{\partial \dot{q}_i \partial \dot{q}_j}(q,\dot{q})v_i v_j = \sum_{i,j = 1}^n m_{ij} v_i v_j > 0
			\end{equation*}
			\noindent for all $(q,\dot{q}) \in TM$ and $v := \sum_{i = 1}^n v_i \partial_i \in T_xM$ such that $v \neq 0$, as $m$ is positive-definite.
		\item The Lagrangian function $L$ is fibrewise supercoersive, that is,
			\begin{equation*}
				\lim_{\norm{\dot{q}}_m \to \infty} \frac{L(q,\dot{q})}{\norm{\dot{q}}_m} = \lim_{\norm{\dot{q}}_m \to \infty} \del[3]{\frac{1}{2}\norm[0]{\dot{q}}_m - \frac{V(q)}{\norm[0]{\dot{q}}_m}} = +\infty
			\end{equation*}
			\noindent for all $q \in M$.
	\end{enumerate}
\end{solution}

\begin{definition}[Symmetry Group]
	A Lie group $G$ is called a \bld{symmetry group of a Lagrangian system $(M,L)$}, if there exists a left action $\psi$ of $G$ on $M$ with 
	\begin{equation*}
		L \circ D\psi_g = L \qquad \forall g \in G.
	\end{equation*}
\end{definition}

\begin{theorem}
	Let $(M,L)$ be a Lagrangian system with symmetry group $G$ and such that the Legendre transform $\tau_L$ is a diffeomorphism. Then $G$ is a symmetry group of the corresponding Hamiltonian system $(T^*M,d\lambda, H_L)$ with
	\begin{equation*}
		\mu(\xi)(q,\tau_L(\dot{q})) = \tau_L(\dot{q})\del[3]{\frac{d}{dt}\bigg\vert_{t = 0}\psi_{\exp(-t\xi)}(q)} \qquad \forall \xi \in \mathfrak{g}, (q,\dot{q}) \in TM,
	\end{equation*}
	\noindent where $\psi$ denotes the smooth left $G$-action on the configuration space $M$. Moreover, the induced action on the phase space $T^*M$ is Hamiltonian and Poisson.
\end{theorem}

\begin{proof}
	Define a smooth left $G$-action $\Psi$ on the phase space $T^*M$ by
	\begin{equation*}
		\Psi_g := \tau_L \circ D\psi_g \circ \tau_L^{-1} \qquad \forall g \in G.
	\end{equation*}
	Applying the Momentum Lemma \ref{lem:momentum_lemma} to this action yields the Theorem. We proceed in five steps.

	\emph{Step 1: $D\psi_g^* \lambda_L = \lambda_L$ for all $g \in G$.} We compute
	\allowdisplaybreaks
	\begin{align*}
		((D\psi_g)^* \lambda_L)(\bar{v}) &= dL((\Phi \circ D\pi_{TM} \circ DD\psi_g)(\bar{v}))\\
		&= dL(\Phi(D\psi_g(\dot{q}),D\pi_{TM} \circ DD\psi_g(\bar{v})))\\
		&= dL(\Phi(D\psi_g(\dot{q}),D\psi_g(D\pi_{TM}(\bar{v}))))\\
		&= dL\del[3]{\frac{d}{d\varepsilon}\bigg\vert_{\varepsilon = 0}D\psi_g(\dot{q} + \varepsilon D\pi_{TM}(\bar{v}))}\\ 
		&= dL \circ D\psi_g\del[3]{\frac{d}{d\varepsilon}\bigg\vert_{\varepsilon = 0}(\dot{q} + \varepsilon D\pi_{TM}(\bar{v}))}\\
		&= dL((\Phi \circ D\pi_{TM})(\bar{v}))\\
		&= \lambda_L(\bar{v})
	\end{align*}
	\noindent for all $\bar{v} \in T_{(q,\dot{q})}TM$.
	
	\emph{Step 2: The induced action $\Psi$ preserves the Liouville form.} For $g \in G$ we compute	
	\allowdisplaybreaks
	\begin{align*}
		\Psi_g^*\lambda &= \del[1]{\tau_L \circ D\psi_g \circ \tau_L^{-1}}^* \lambda\\
		&= \del[1]{\tau_L^{-1}}^* (D\psi_g)^* \tau_L^* \lambda\\
		&= \del[1]{\tau_L^{-1}}^* (D\psi_g)^* \lambda_L\\
		&= \del[1]{\tau_L^{-1}}^* \lambda_L\\ 
		&= \lambda
	\end{align*}
	\noindent by Step 1.

	\emph{Step 3: The momentum map is of the stated form.} We compute
	\allowdisplaybreaks
	\begin{align*}
		\mu(\xi)(q,\tau_L(\dot{q})) &= i_{\hat{\xi}}\lambda(q,\tau_L(\dot{q}))\\
		&= \lambda_{(q,\tau_L(\dot{q}))}\del[3]{\frac{d}{dt}\bigg\vert_{t = 0}\Psi_{\exp(-t\xi)}(q,\tau_L(\dot{q}))}\\
		&= p\del[3]{\frac{d}{dt}\bigg\vert_{t = 0} \pi_{T^*M} \circ \Psi_{\exp(-t\xi)}(q,\tau_L(\dot{q}))}\\
		&= p\del[3]{\frac{d}{dt}\bigg\vert_{t = 0} \pi_{TM} \circ D\psi_{\exp(-t\xi)} \circ \tau^{-1}_L(q,\tau_L(\dot{q}))}\\
		&= p\del[3]{\frac{d}{dt}\bigg\vert_{t = 0} \psi_{\exp(-t\xi)} \circ \pi_{TM} \circ \tau^{-1}_L(q,\tau_L(\dot{q}))}\\
		&= p\del[3]{\frac{d}{dt}\bigg\vert_{t = 0} \psi_{\exp(-t\xi)} \circ \pi_{T^*M}(q,\tau_L(\dot{q}))}\\
		&= \tau_L(\dot{q})\del[3]{\frac{d}{dt}\bigg\vert_{t = 0} \psi_{\exp(-t\xi)}(q)}
	\end{align*}
	\noindent for all $\xi \in \mathfrak{g}$ and $(q,\dot{q}) \in TM$.

	\emph{Step 4: $E_L \circ D\psi_g = E_L$ for all $g \in G$.} We compute
		\begin{align*}
		E_L(\psi_g(q),D\theta_g(\dot{q})) &= \tau_L(D\psi_g(\dot{q})(D\psi_g(\dot{q})) - L \circ D\psi_g(q,\dot{q})\\
		&= \frac{d}{d\varepsilon}\bigg \vert_{\varepsilon = 0} L_{\psi_g(q)} \circ D\psi_g((1 + \varepsilon)\dot{q}) - L(q,\dot{q})\\
		&= \frac{d}{d\varepsilon}\bigg \vert_{\varepsilon = 0} L_q((1 + \varepsilon)\dot{q}) - L(q,\dot{q})\\
		&= E_L(q,\dot{q})
	\end{align*}
	\noindent for all $g \in G$ and $(q,\dot{q}) \in TM$.

	\emph{Step 5: $H_L \circ \Psi_g = H_L$ for all $g \in G$.} Using Step 4 we conclude
	\allowdisplaybreaks
	\begin{equation*}
		H_L \circ \Psi_g = E_L \circ D\psi_g \circ \tau_L^{-1} = E_L \circ \tau_L^{-1} = H_L
	\end{equation*}
	\noindent for all $g \in G$. 
\end{proof}

\subsection{The Lagrangian and Hamiltonian Action Functional}

\begin{definition}[Lagrangian Action Functional]
For a Lagrangian system $(M,L)$ the corresponding \bld{Lagrangian action functional} is
\begin{equation*}
    \mathscr{E}_L \colon \mathscr{L}M \to \mathbb{R}, \qquad \mathscr{E}_L(\gamma) := \int_0^1 L(\dot{\gamma}(t))dt,
\end{equation*}
\noindent where $\mathscr{L}M := C^\infty(\mathbb{T},M)$ denotes the loop space of $M$ with $\mathbb{T} = \mathbb{R} /\mathbb{Z}$.
\end{definition}

\begin{exercise}
    Let $(M,L)$ a Lagrangian system such that the Legendre transform $\tau_L$ is a diffeomorphism. Show that
    \begin{equation*}
        \mathscr{E}_L(\gamma) = \int_0^1 (\tau_L \circ \dot{\gamma})^*\lambda - \int_0^1 H_L(\tau_L(\dot{\gamma}(t)))dt
    \end{equation*}
    \noindent for all $\gamma \in \mathscr{L}M$.
\end{exercise}

\begin{solution}
    We compute using the definitions
    \allowdisplaybreaks
    \begin{align*}
        \mathscr{E}_L(\gamma) &= \int_0^1 L(\dot{\gamma}(t))dt\\
        &= \int_0^1 \tau_L(\dot{\gamma})(\dot{\gamma}) - \int_0^1 E_L(\dot{\gamma}(t))dt\\
        &= \int_0^1 \tau_L(\dot{\gamma})(\dot{\gamma}) - \int_0^1 H_L(\tau_L(\dot{\gamma}(t)))dt\\
        &= \int_0^1 \dot{\gamma}^*\lambda_L - \int_0^1 H_L(\tau_L(\dot{\gamma}(t)))dt\\
        &= \int_0^1 (\tau_L \circ \dot{\gamma})^*(\tau_L)_*\lambda_L - \int_0^1 H_L(\tau_L(\dot{\gamma}(t)))dt\\
        &= \int_0^1 (\tau_L \circ \dot{\gamma})^*\lambda - \int_0^1 H_L(\tau_L(\dot{\gamma}(t)))dt.
    \end{align*}
\end{solution}

\begin{definition}[Hamiltonian Action Functional]
    Let $(M,\lambda,H)$ be an exact Hamiltonian system. Then the corresponding \bld{Hamiltonian action functional} is
\begin{equation*}
    \mathscr{A}_H \colon \mathscr{L}M \to \mathbb{R}, \qquad \mathscr{A}_H(\gamma) := \int_0^1 \gamma^*\lambda - \int_0^1 H(\gamma(t))dt.
\end{equation*}
\end{definition}

The following principle underlies classical mechanics \cite[Chapter~8]{knauf:cm:2018}.

\begin{axiom}[Hamilton's Principle of Extremal Action]
    \label{ax:hamilton}
    A path $\gamma \in \mathscr{L}M$ describes a motion of a Lagrangian system $(M,L)$, if
    \begin{equation*}
        \frac{d}{d\varepsilon}\bigg\vert_{\varepsilon = 0}\mathscr{E}_L(\gamma_\varepsilon) = 0
    \end{equation*}
    \noindent for all variations $\gamma_\varepsilon$ of $\gamma$, that is, given $X \in \mathfrak{X}(\gamma)$, we set
    \begin{equation*}
        \gamma_\varepsilon(t) := \exp^\nabla_{\gamma(t)}(\varepsilon X(t)) \qquad \forall t \in \mathbb{T}
    \end{equation*}
    \noindent for some Levi--Civita connection $\nabla$.
\end{axiom}

\begin{remark}
    Axiom \ref{ax:hamilton} is just the statement that motions of Lagrangian systems are critical points $\gamma \in \operatorname{Crit}(\mathscr{E}_L)$, where the differential $d\mathscr{E}_L$ of the Lagrangian action functional $\mathscr{E}_L$ is defined by
    \begin{equation}
        d\mathscr{E}_L \colon \mathscr{L}M \to T^*\mathscr{L}M, \qquad d\mathscr{E}_L(\gamma)(X) =  \frac{d}{d\varepsilon}\bigg\vert_{\varepsilon = 0}\mathscr{E}_L(\gamma_\varepsilon),
    \end{equation}
    \noindent with infinitesimal variation
    \begin{equation*}
        \frac{d}{d\varepsilon}\bigg\vert_{\varepsilon = 0}\gamma_\varepsilon = X.
    \end{equation*}
\end{remark}

\begin{lemma}
    \label{lem:el}
    In standard coordinates $(q,\dot{q})$ on the tangent bundle $TM$, the differential $d\mathscr{E}_L$ of the Lagrangian action functional $\mathscr{E}_L$ is given by
    \begin{equation*}
        d\mathscr{E}_L(\gamma)(X) = \int_0^1 \del[3]{\frac{\partial L}{\partial q} - \frac{d}{dt}\frac{\partial L}{\partial \dot{q}}}X(t)dt.
    \end{equation*}
    In particular, $d\mathscr{E}_L$ is well-defined and any motion of a Lagrangian system satisfies locally the \bld{Euler--Lagrange equations}
    \begin{equation}
        \label{eq:el}
        \frac{d}{dt}\frac{\partial L}{\partial \dot{q}} = \frac{\partial L}{\partial q}.
    \end{equation}
\end{lemma}

\begin{proof}
We compute
\allowdisplaybreaks
\begin{align*}
    d\mathscr{E}_L(\gamma)(X) &= \frac{d}{d\varepsilon}\bigg\vert_{\varepsilon = 0} \mathscr{E}_L(\gamma_\varepsilon)\\
    &= \int_0^1 \frac{d}{d\varepsilon}\bigg\vert_{\varepsilon = 0} L(\dot{\gamma}_\varepsilon(t))dt\\
    &= \int_0^1 dL\del[3]{\frac{d}{d\varepsilon}\bigg\vert_{\varepsilon = 0} \dot{\gamma}_\varepsilon(t)}dt\\
    &= \int_0^1 \del[3]{\frac{\partial L}{\partial q} \frac{d\gamma_\varepsilon(t)}{d \varepsilon} \bigg\vert_{\varepsilon = 0} + \frac{\partial L}{\partial \dot{q}} \frac{d\dot{\gamma}_\varepsilon(t)}{d \varepsilon} \bigg\vert_{\varepsilon = 0}}dt\\
    &= \int_0^1 \del[3]{\frac{\partial L}{\partial q} \frac{d\gamma_\varepsilon(t)}{d \varepsilon} \bigg\vert_{\varepsilon = 0} + \frac{\partial L}{\partial \dot{q}}\frac{d}{dt} \frac{d\gamma_\varepsilon(t)}{d \varepsilon} \bigg\vert_{\varepsilon = 0}}dt\\
    &= \int_0^1 \del[3]{\frac{\partial L}{\partial q}X(t) + \frac{\partial L}{\partial \dot{q}}\frac{d}{dt}X(t)}dt\\
    &= \int_0^1 \frac{\partial L}{\partial q}X(t)dt + \int_0^1 \frac{\partial L}{\partial \dot{q}}\frac{d}{dt}X(t)dt\\
    &= \int_0^1 \frac{\partial L}{\partial q}X(t)dt + \frac{\partial L}{\partial \dot{q}}X(t)\big\vert_0^1 - \int_0^1 \frac{d}{dt}\frac{\partial L}{\partial \dot{q}}X(t)dt\\
    &= \int_0^1 \frac{\partial L}{\partial q}X(t)dt - \int_0^1 \frac{d}{dt}\frac{\partial L}{\partial \dot{q}}X(t)dt\\
    &= \int_0^1 \del[3]{\frac{\partial L}{\partial q} - \frac{d}{dt}\frac{\partial L}{\partial \dot{q}}}X(t)dt
\end{align*}
\noindent using integration by parts. As $X \in \mathfrak{X}(\gamma)$ is arbitrary, we conclude \eqref{eq:el} by Exercise \ref{ex:fundamental}.
\end{proof}

\begin{exercise}[Fundamental Lemma of Calculus of Variations]
    \label{ex:fundamental}
    Let $M$ be a smooth manifold and $f \in C^\infty(M)$. Prove that, if
    \begin{equation*}
        \int_M f\varphi = 0 \qquad \forall \varphi \in C^\infty_c(M),
    \end{equation*}
    \noindent then $f = 0$.
\end{exercise}

\begin{solution}
    Assume that $f \neq 0$. Then there exists $x_0 \in M$ such that $f(x_0) \neq 0$. Without loss of generality, we may assume that $f(x_0) > 0$. By continuity of $f$, there exists $\delta > 0$ with
    \begin{equation*}
        f(x) \in B_{f(x_0)/2}(f(x_0)) \qquad \forall x \in B_\delta(x_0).
    \end{equation*}
    Choose a smooth bump function $\varphi \in C_c^\infty(M)$ supported in $B_\delta(x_0)$ and with $\varphi = 1$ on $\bar{B}_{\delta/2}(x_0)$. Then
    \begin{equation*}
        0 = \int_M f\varphi = \int_{B_\delta(x_0)}f\varphi \geq \int_{B_{\delta/2}(x_0)} f\varphi > \frac{1}{2}f(x_0)\operatorname{Vol}(B_{\delta/2}(x_0)) > 0.
    \end{equation*}
\end{solution}

A more elegant way of writing the Euler--Lagrange equations \eqref{eq:el} is via the Hamiltonian formalism.

\begin{lemma}
    \label{lem:hamiltonian_equations}
    A path $\gamma \in \mathscr{L}M$ is a motion of an exact Hamiltonian system $(M,\lambda,H)$ if and only if it is an integral curve of the Hamiltonian vector field $X_H$, that is, $\gamma$ satisfies the ordinary differential equation
    \begin{equation}
        \label{eq:ham}
        \dot{\gamma}(t) = X_H(\gamma(t)) \qquad \forall t \in \mathbb{T}.
    \end{equation}
\end{lemma}

\begin{proof}
    We compute
    \allowdisplaybreaks
    \begin{align*}
        d\mathscr{A}_H(\gamma)(X) &= \frac{d}{d\varepsilon}\bigg\vert_{\varepsilon = 0} \mathscr{A}_H(\gamma_\varepsilon)\\
        &= \int_0^1 \frac{d}{d\varepsilon}\bigg\vert_{\varepsilon = 0} \gamma_\varepsilon^*\lambda - \int_0^1 \frac{d}{d\varepsilon}\bigg\vert_{\varepsilon = 0} H(\gamma_\varepsilon(t))dt\\
        &= \int_0^1 \gamma^*L_X\lambda - \int_0^1 dH(X(t))dt\\
        &= \int_0^1 \gamma^*(i_Xd\lambda + di_X\lambda) - \int_0^1 dH(X(t))dt\\
        &= \int_0^1 d\lambda(X(t),\dot{\gamma}(t))dt + \lambda(X)\big\vert_0^1 - \int_0^1 dH(X(t))dt\\
        &= \int_0^1 d\lambda(X(t),\dot{\gamma}(t))dt + \int_0^1 d\lambda(X_H(\gamma(t)),X(t))dt\\
        &= \int_0^1 d\lambda(X(t),\dot{\gamma}(t) - X_H(\gamma(t)))dt
    \end{align*}
    \noindent using Fisherman's formula \ref{lem:fisherman} and Stoke's theorem. As $X \in \mathfrak{X}(\gamma)$ is arbitrary, we conclude \eqref{eq:ham} by Exercise \ref{ex:fundamental} and the nondegeneracy of the symplectic form $d\lambda$.
\end{proof}

\begin{exercise}
    Compare the proof of Lemma \ref{lem:el} with the proof of Lemma \ref{lem:hamiltonian_equations}. What are similarities and differences? Which proof is more efficient?
\end{exercise}

\begin{solution}
    Both proofs make use of the G\^ateaux derivative of a functional defined on an infinite-dimensional space. There is no need to consider the tangent space of the infinite-dimensional manifold $\mathscr{L}M$ in both cases, as one can elegantly define the tangent space at a loop $\gamma$ as vector fields along $\gamma$. As $\mathscr{L}M$ is only a Fr\'echet manifold and not even a Banach manifold, this would introduce unnecessary difficulties. In Floer theory one usually works with the Hilbert manifold $W^{k,2}(\mathbb{T},M)$ anyway. The proof of Lemma \ref{lem:el} is local in nature, appealing to coordinates, whereas the proof of Lemma \ref{lem:hamiltonian_equations} is global. Both perspectives are important in mathematics and physics, as one should be able to grasp the over spanning concept as well as being able to perform dirty computations in coordinates. Thus both proofs are efficient in some sense.
\end{solution}

\begin{exercise}
    Give a proof of Lemma \ref{lem:hamiltonian_equations} in the case $(T^*M,p dq,H_L)$ for a Lagrangian function $L \in C^\infty(TM)$ for a configuration space $M$. 
\end{exercise}

\begin{solution}
    In coordinates we have that
	\begin{equation*}
		H_L(q,p) = \del[3]{\frac{\partial L}{\partial \dot{q}}\dot{q} - L(q,\dot{q})}\bigg\vert_{p = \frac{\partial L}{\partial \dot{q}}}.
	\end{equation*}
	Therefore
	\begin{equation*}
		\frac{\partial H_L}{\partial p} = \frac{\partial}{\partial p} \del[1]{p \dot{q} - L(q,\dot{q})}\big\vert_{p = \frac{\partial L}{\partial \dot{q}}} = \dot{q}
	\end{equation*}
	\noindent and
    \begin{equation*}
		\frac{\partial H_L}{\partial q} = \frac{\partial}{\partial q}\del[3]{p\dot{q} - L(q,\dot{q})}\bigg\vert_{p = \frac{\partial L}{\partial \dot{q}}} = - \frac{\partial L}{\partial q}\bigg\vert_{p = \frac{\partial L}{\partial \dot{q}}} = -\frac{d}{dt}\frac{\partial L}{\partial \dot{q}} \bigg\vert_{p = \frac{\partial L}{\partial \dot{q}}} = -\dot{p}
	\end{equation*}
    \noindent using the Euler--Lagrange equations \eqref{eq:el}.
\end{solution}

\begin{exercise}
    Pair up and compare and contrast the Hamiltonian formalism to the Lagragian formalism by filling out the table below.
\end{exercise}

\begin{summary*}{The Hamiltonian and the Lagrangian Formalism}.
    \renewcommand{\arraystretch}{1.5}
    \begin{tabular}{c|c|c}
        & \textbf{Hamiltonian Formalism} & \textbf{Lagrangian Formalism}\\
        \hline
        \emph{Configuration} & Smooth manifold & Smooth manifold\\
        \emph{space} & $M^n$ & $M^n$\\
        \hline
        \emph{Function} & Hamiltonian function  & Larangian function \\
        & $H \in C^\infty(T^*M)$ & $L \in C^\infty(TM)$\\
        \hline
        \emph{Domain} & Phase space cotangent  & State space tangent \\
        & bundle $T^*M$   & bundle $TM$\\
        \hline
        \emph{Coordinates} & $(q,p)$ for $q$ position & $(q,\dot{q})$ for $q$ position\\
        & and $p$ momentum & and $\dot{q}$ velocity\\
        \hline
        \emph{Dimension} & $2n$ & $2n$\\
        \hline
        \emph{Equations of} & Hamiltonian equations & Euler--Lagrange Equations\\
        \emph{motion} & $\displaystyle\dot{q} = \frac{\partial H}{\partial p}$ and $\displaystyle\dot{p} = -\frac{\partial H}{\partial q}$ & $\displaystyle \frac{d}{dt}\frac{\partial L}{\partial\dot{q}} = \frac{\partial L}{\partial q}$\\
        \hline
        \emph{Functional} & Hamiltonian action $\mathscr{A}_H$ & Lagrangian action $\mathscr{E}_L$\\
        \hline
        \emph{Hamilton's} & $(q,p) \in \operatorname{Crit}(\mathscr{A}_H)$, that is & $q \in \Crit(\mathscr{E}_L)$, that is\\
        \emph{principle} & $\displaystyle \frac{d}{d\varepsilon}\bigg\vert_{\varepsilon = 0}\mathscr{A}_H(q_\varepsilon,p_\varepsilon) = 0$ & $\displaystyle\frac{d}{d\varepsilon}\bigg\vert_{\varepsilon = 0}\mathscr{E}_L(q_\varepsilon) = 0$\\
        & for all variations $(q_\varepsilon,p_\varepsilon)$ & for all variations $q_\varepsilon$\\
        \hline
        \emph{Transformation} & Inverse Legendre transformation & Legendre transformation\\
        & $\tau_L^{-1} \colon T^* M \to TM$ & $\tau_L \colon TM \to T^*M$\\
        \hline
        \emph{Symmetries} & Symmetry group $G$ with & Symmetry group $G$ with\\
        & $\forall g \in G : H \circ \psi_g = H$ & $\forall g \in G : L \circ D\psi_g = L$\\
        \hline
        \emph{Momentum} & $\mu \colon \mathfrak{g} \to C^\infty(T^*M)$ & $\mu \circ \tau_L^* \colon \mathfrak{g} \to C^\infty(TM)$,\\
        \emph{map} & & $\tau_L^* \colon C^\infty(T^*M) \to C^\infty(TM)$\\
        \hline
    \end{tabular}                
\end{summary*}
\newpage
\section{Regular Energy Surfaces}
\begin{question}
    Let $(M,\omega,H)$ be a Hamiltonian system. Assume that $X_H$ is complete, meaning that the flow
    \begin{equation*}
        \phi \colon \mathbb{R} \times M \to M, \qquad (t,x) \to \phi_t(x),
    \end{equation*}
    \noindent defined by
    \begin{equation*}
        \phi_0 = \id_M \qquad \text{and} \qquad \frac{d}{dt}\phi_t = X_H \circ \phi_t
    \end{equation*}
    \noindent is global. Which options are correct?
    \begin{itemize}
        \item[$\square$] The energy $H$ is constant along flow lines, meaning that
        \begin{equation*}
            \frac{d}{dt}H \circ \phi_t = 0 \qquad \forall t \in \mathbb{R}.
        \end{equation*}
        \item[$\square$] Let $\varphi \in \Symp(M,\omega)$ be a symplectomorphism, meaning that $\varphi \in \Diff(M)$ is a diffeomorphism and $\varphi^*\omega = \omega$. Then we have that
        \begin{equation*}
            \varphi \circ \phi_t^{X_{\varphi^*H}} = \phi_t^{X_H} \circ \varphi.
        \end{equation*}
        \item[$\square$] The map
        \begin{equation*}
            \del[1]{C^\infty(M),\{ \cdot,\cdot\}} \to \del[1]{\mathfrak{X}(M),[\cdot,\cdot]}, \qquad f \mapsto X_f
        \end{equation*}
        \noindent is a Lie algebra homomorphism, where
        \begin{equation*}
            \{\cdot,\cdot\}\colon C^\infty(M) \times C^\infty(M) \to C^\infty(M), \qquad \{f,g\} := \omega(f,g)
        \end{equation*}
        \noindent denotes the Poisson bracket.
        \item[$\square$] Let $\psi \colon G \times M \to M$ be a smooth action by a Lie group $G$ with corresponding Lie algebra $\mathfrak{g}$. Assume that the map
        \begin{equation*}
            \mu \colon \mathfrak{g} \to C^\infty(M), \qquad \mu(\xi) := \frac{d}{dt}\bigg\vert_{t = 0}\psi_{\exp(-t\xi)}
        \end{equation*}
        \noindent is linear. Then we have that
        \begin{equation*}
            \{\mu(\xi),H\} = 0 \qquad \forall \xi \in \mathfrak{g}.
        \end{equation*}
    \end{itemize}
\end{question}

\begin{solution}
    Correct is the first and the second option.
\end{solution}

\subsection{The Poincar\'e Recurrence Theorem}

\begin{definition}[Regular Energy Surface]
	A \bld{regular energy surface in a Hamiltonian system $(M,\omega,H)$} is defined to be an embedded hypersurface $\Sigma = H^{-1}(0)$ such that $\Crit(H) \cap \Sigma = \emptyset$.
\end{definition}

On compact regular energy surfaces, there do exist special finite regular Borel measures \cite[Proposition~V.31]{zehnder:ds:2010}. Let $\Sigma$ be a compact regular energy surface in a Hamiltonian system $(M^{2n},\omega,H)$ and $\varphi \in \Symp(M,\omega)$ such that $H \circ \varphi = H$. Then we can turn $\varphi\vert_\Sigma\colon \Sigma \mapsto \Sigma$ into a measure-theoretical discrete dynamical system. In order to do that, we need to construct a suitable $\varphi\vert_\Sigma$-invariant probability measure on $\Sigma$. This measure should be well-behaved, in the sense that it comes from an $\omega$-induced volume form on $\Sigma$. If $\Omega_\Sigma \in \Omega^{2n - 1}(\Sigma)$ is a volume form on $\Sigma$ such that $\varphi\vert_\Sigma^*\Omega_\Sigma = \Omega_\Sigma$, there exists a unique regular $\varphi\vert_\Sigma$-invariant Borel probability measure $\mu_\Sigma$ on $\Sigma$, such that
\begin{equation*}
	\int_\Sigma f\bar{\Omega}_\Sigma = \int_M f\mu_\Sigma \qquad \textstyle\bar{\Omega}_\Sigma := \Omega_\Sigma/\int_\Sigma \Omega_\Sigma,
\end{equation*}
\noindent for all $f \in C(\Sigma)$. Indeed, this immediately follows from the Stone--Weierstrass Theorem \cite[p.~392]{cohn:mt:2013} and the Riesz Representation Theorem \cite[Theorem~7.2.8]{cohn:mt:2013}.

Let us construct a $\varphi\vert_\Sigma$-invariant volume form $\Omega_\Sigma$ on $\Sigma$. As $dH \neq 0$ on $\Sigma$ by definition of a regular energy surface, we find a neighbourhood $U$ of $\Sigma$ in $M$ such that $dH \neq 0$ on $U$ as $\Sigma$ is assumed to be compact. Pick $\alpha \in \Omega^{2n - 1}(U)$ such that
\begin{equation}
	\omega^n = dH \wedge \alpha.
	\label{eq:volume_form}
\end{equation}
Such a form $\alpha$ does exist, as one may take $\alpha = i_X\omega^n$, where
\begin{equation*}
	X := \grad H /\norm[0]{\grad H}^2
\end{equation*}
\noindent with respect to some Riemannian metric on $U$. Then the volume form $\iota_\Sigma^* \alpha$ is uniquely determined by the requirement \eqref{eq:volume_form}. Indeed, if $\alpha' \in \Omega^{2n - 1}(U)$ satisfies \eqref{eq:volume_form}, there exists a differential form $\beta \in \Omega^{2n - 2}(U)$ such that
\begin{equation}
	\alpha - \alpha' = dH \wedge \beta.
	\label{eq:volume_form_induced}
\end{equation}
For example, take $i_X(\alpha - \alpha')$. Applying $\iota_\Sigma^*$ to \eqref{eq:volume_form_induced} yields $\iota_\Sigma^*\alpha = \iota_\Sigma^*\alpha'$. Let us denote by $\Omega_\Sigma := \iota_\Sigma' \alpha$, where $\alpha \in \Omega^{2n - 1}(U)$ is any differential form satisfying \eqref{eq:volume_form}. We compute
\begin{equation*}
	\omega^n = (\varphi^*\omega)^n = \varphi^*\omega^n = d(H \circ \varphi) \wedge \varphi^*\alpha = dH \wedge \varphi^*\alpha.
\end{equation*}
Consequently, by uniqueness 
\begin{equation*}
	\Omega_\Sigma = \iota_\Sigma^*\alpha = \iota_\Sigma^*\varphi^*\alpha = (\varphi \circ \iota_\Sigma)^*\alpha = (\iota_\Sigma \circ \varphi\vert_\Sigma)^*\alpha = \varphi\vert_\Sigma^*\Omega_\Sigma,
\end{equation*}
\noindent and thus $\varphi\vert_\Sigma \colon \Sigma \to \Sigma$ is a measure-theoretical dynamical system on $(\Sigma,\mathcal{B}(\Sigma),\mu_\Sigma)$.

\begin{lemma}[{Poincar\'e's Recurrence Theorem, \cite[I.15]{zehnder:ds:2010}}]	
	\label{lem:poincare_recurrence_theorem}
	Let $\Sigma$ be a compact regular energy surface in a Hamiltonian system and $\varphi \in \Symp(M,\omega)$ such that $H \circ \varphi = H$. Then for $\mu_\Sigma$-almost every point $x \in \Sigma$ there exists a sequence $(k_\nu) \subseteq \mathbb{N}$ such that
	\begin{equation*}
		k_\nu \to \infty \qquad \text{and} \qquad \lim_{\nu \to \infty} \varphi^{k_\nu}(x) = x.
	\end{equation*}
\end{lemma}

\begin{proof}
	A routine computation shows 
	\begin{equation}
		\label{eq:measure_recurrence}
		\mu_\Sigma\del[4]{A \cap \bigcap_{k \geq 0}\bigcup_{l \geq k}\varphi^{-l}(A)} = \mu_\Sigma(A) \qquad \forall A \in \mathcal{B}(\Sigma).
	\end{equation}
	Fix a Riemannian metric on $\Sigma$. Then $\Sigma$ is a compact metric space with respect to the induced Riemannian distance function and metric topology coinciding with the manifold topology. For every $n \in \mathbb{N}$ there exists a finite index set $I_n$ such that $\del[1]{B_{1/n}(x_{i,n})}_{i \in I_n}$ is an open cover for $\Sigma$. Define 
	\begin{equation*}
		N := \bigcup_{n \in \mathbb{N}}\bigcup_{i \in I_n} \del[4]{B_{1/n}(x_{i,n}) \setminus \del[4]{B_{1/n}(x_{i,n}) \cap \bigcap_{k \geq 0}\bigcup_{l \geq k} \varphi^{-l}\del[1]{B_{1/n}(x_{i,n})}}}.
	\end{equation*}
	Then $N \in \mathcal{B}(\Sigma)$ and $\mu(N) = 0$. Indeed, we have that
	\begin{align*}
		\mu(N) &\leq \sum_{n \in \mathbb{N}} \sum_{i \in I_n} \mu\del[4]{B_{1/n}(x_{i,n}) \setminus \del[4]{B_{1/n}(x_{i,n}) \cap \bigcap_{k \geq 0}\bigcup_{l \geq k} \varphi^{-l}\del[1]{B_{1/n}(x_{i,n})}}}\\
		&= \sum_{n \in \mathbb{N}} \sum_{i \in I_n} \mu\del[1]{B_{1/n}(x_{i,n})} - \mu\del[4]{B_{1/n}(x_{i,n}) \cap \bigcap_{k \geq 0}\bigcup_{l \geq k} \varphi^{-l}\del[1]{B_{1/n}(x_{i,n})}}\\
		&= 0, 
	\end{align*}
	\noindent by \eqref{eq:measure_recurrence}. Moreover, every $x \in N^c$ satisfies the condition of the theorem. Indeed, $x \in N^c$ means that for all $n \in \mathbb{N}$ and $i \in I_n$
	\begin{equation*}
		x \in \del[1]{B_{1/n}(x_{i,n})}^c \qquad \text{or} \qquad x \in \bigcap_{k \geq 0}\bigcup_{l \geq k} \varphi^{-l}\del[1]{B_{1/n}(x_{i,n})}.
	\end{equation*}
	Since $\del[1]{B_{1/n}(x_{i,n})}_{i \in I_n}$ is an open cover for $S_c$, we conclude that 
	\begin{equation}
		\label{eq:contained_in}
		x \in \bigcap_{k \geq 0}\bigcup_{l \geq k} \varphi^{-l}\del[1]{B_{1/n}(x_{i,n})}
	\end{equation}
	\noindent for some $i \in I_n$. Consequently, for every $n \in \mathbb{N}$ there exists an index $i_n \in I_n$ such that \eqref{eq:contained_in} holds.
\end{proof}

Using preservation of energy \ref{cor:preservation_of_energy} yields the following corollary.

\begin{corollary}
    \label{cor:flow}
	Let $\Sigma$ be a compact regular energy surface in a Hamiltonian system $(M,\omega,H)$. Denote by $\varphi := \phi_1^{X_H} \in \Symp(M,\omega)$ the flow of $X_H$ at time $t = 1$. For a Borel measurable set $A \in \mathcal{B}(\Sigma)$ with $\mu_\Sigma(A) > 0$ we define the \bld{Poincar\'e return map}
    \begin{equation}
        \label{eq:return_map}
        P_A \colon A \to A, \qquad P_A(x) := \varphi^{\tau_A(x)}(x),
    \end{equation}
    \noindent where
    \begin{equation*}
        \tau_A \colon A \to \mathbb{N} \cup \{\infty\},\qquad \tau_A(x) := \inf\cbr[1]{k \geq 1 : \varphi^k(x) \in A} 
    \end{equation*}
    \noindent is the \bld{first return time to $A$}. Then the Poincar\'e return map $P_A$ is well-defined $\mu_\Sigma$-almost everywhere.
\end{corollary}

\begin{exercise}
    Prove Corollary \ref{cor:flow} and sketch the situation of the statement. In particular, explain why the Poincar\'e recurrence theorem \ref{lem:poincare_recurrence_theorem} is called recurrence theorem!
\end{exercise}

\begin{figure}[h!tb]
    \centering
    \includegraphics[width=.8\textwidth]{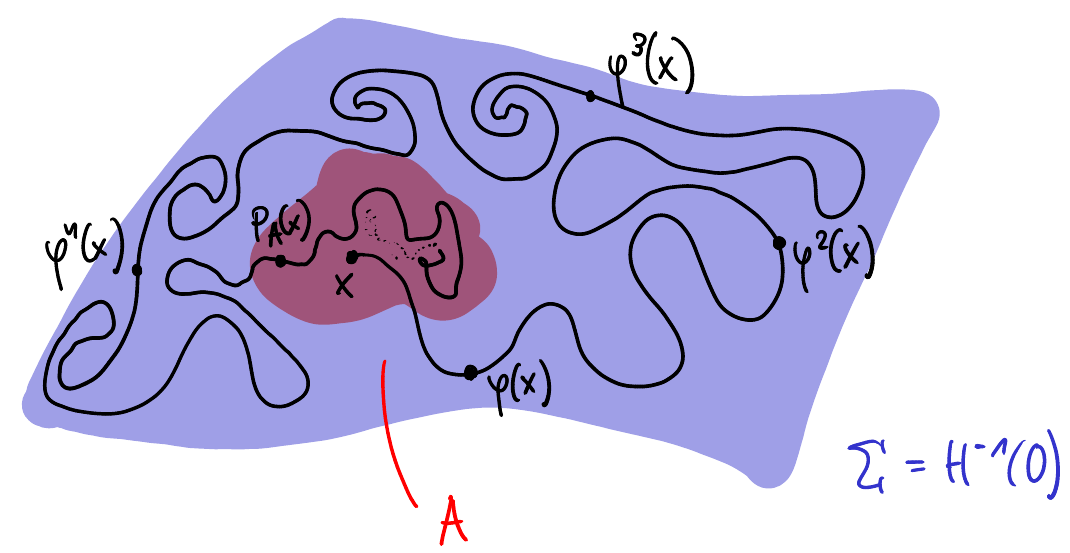}
    \caption{An illustration of the recurrence property of measure preserving dynamical systems}
\end{figure}

\begin{solution}
    By preservation of energy \ref{cor:preservation_of_energy}, we have that $H \circ \varphi = H$. As $\mathcal{B}(\Sigma)$ is generated by open sets, we can without loss of generality assume that $A$ is open. Pick $x \in A$. Thus there exists a ball $B_\varepsilon(x)$ contained in $A$. By the Poincar\'e recurrence theorem \ref{lem:poincare_recurrence_theorem}, $\mu_\Sigma$-almost every point in $B_\varepsilon(x)$ returns to $B_\varepsilon(x)$ after finitely many iterations of $\varphi$. Thus $\tau_A(x)$ is finite for almost every point in $B_\varepsilon(x)$ and thus also for $A$.
\end{solution}

\begin{remark}[Poincar\'e's Continuation Method]
    The Poincar\'e return map \eqref{eq:return_map} is important in the study of the existence of periodic orbits on compact energy hypersurfaces in Hamiltonian systems. For a more detailed introduction see \cite[Section~VIII.4]{zehnder:ds:2010} and \cite[Section~8.2]{abrahammarsden:cm:1978}.
\end{remark}

\subsection{Stable Hypersurfaces}

\begin{definition}[{Hamiltonian Manifold, \cite[Definition~2.4.1]{frauenfelderkoert:3bp:2018}}]
	A \bld{Hamiltonian manifold}\index{Manifold!Hamiltonian} is defined to be a pair $(\Sigma,\omega)$, where $\Sigma$ is an odd-dimensional smooth manifold and $\omega \in \Omega^2(\Sigma)$ is closed such that $\ker \omega$ is a line distribution. The foliation inducing the line distribution $\ker \omega$ is called the \bld{characteristic foliation}.
\end{definition}

\begin{exercise}[Regular Energy Surface]
    \label{ex:hamiltonian_manifold}
	Let $\Sigma$ be a regular energy surface in a Hamiltonian system $(M,\omega,H)$. Show that $(\Sigma,\omega\vert_\Sigma)$ is a Hamiltonian manifold and that $\ker \omega\vert_\Sigma$ is spanned by the Hamiltonian vector field $X_H\vert_\Sigma \colon \Sigma \to T\Sigma$.
\end{exercise}

\begin{solution}
    The hypersurface $\Sigma^{2n - 1} \subseteq M^{2n}$ is clearly of odd dimension and $\omega\vert_\Sigma$ is closed. It remains to show that $\ker \omega\vert_\Sigma \to \Sigma$ is a line bundle spanned by $X_H\vert_\Sigma$. Using the dimension formula \cite[Lemma~2.1.1]{mcduffsalamon:st:2017}
    \begin{equation*}
        \dim T_x \Sigma + \dim (T_x \Sigma)^\omega = \dim T_x M = 2n
    \end{equation*}
    \noindent yields $\dim (T_x \Sigma)^\omega = 1$ for the symplectic complement
    \begin{equation*}
    (T_x \Sigma)^\omega = \cbr[1]{v \in T_x M : i_v\omega\vert_{T_x \Sigma} = 0}. 
    \end{equation*}
    Now we have that $X_H(x) \in (T_x \Sigma)^\omega$ for every $x \in \Sigma$ as 
    \begin{equation*}
        \omega(X_H(x),v) = -dH_x(v) = 0 \qquad \forall v \in T_x\Sigma = \ker dH_x.
    \end{equation*}
    Also $X_H(x) \in T_x\Sigma$ as
    \begin{equation*}
        dH_x(X_H(x)) = \omega(X_H(x),X_H(x)) = 0.
    \end{equation*}
    Consequently 
    \begin{equation*}
        \ker \omega\vert_\Sigma = (T\Sigma)^\omega = \langle X_H\vert_\Sigma \rangle.
    \end{equation*}
\end{solution}

\begin{definition}[{Stable Hamiltonian Manifold, \cite[p.~1773]{cieliebakfrauenfelderpaternain:mane:2010}}]
    \label{def:stable}
	A Hamiltonian manifold $(\Sigma,\omega)$ is called \bld{stable}\index{Manifold!Hamiltonian!stable}, if there exists $\lambda \in \Omega^1(\Sigma)$ which is nowhere-vanishing on $\ker \omega$ and such that $\ker \omega \subseteq \ker d\lambda$. We write $(\Sigma,\omega,\lambda)$ for a stable Hamiltonian manifold.
\end{definition}

\begin{exercise}[Regular Energy Hypersurface]
    \label{ex:stable}
	Let $\Sigma$ be a regular energy surface in a Hamiltonian system $(M,\omega,H)$ such that there exists a vector field $X$ in a neighbourhood of $\Sigma$ with $X$ being nowhere tangent to $\Sigma$ and $\ker \omega\vert_\Sigma \subseteq \ker L_X\omega \vert_\Sigma$. Show that $(\Sigma,\omega\vert_\Sigma,i_X\omega\vert_\Sigma)$ is a stable Hamiltonian manifold.
\end{exercise}

\begin{solution}
    From Exercise \ref{ex:hamiltonian_manifold} if follows that $(\Sigma,\omega\vert_\Sigma)$ is a Hamiltonian manifold. The condition that $\lambda \in \Omega^1(\Sigma^{2n - 1})$ is a stabilising form, is equivalent to the fact that $\lambda \wedge \omega^{n - 1} \in \Omega^{2n - 1}(\Sigma)$ is a volume form and thus nowhere vanishing on $\Sigma$. For $\lambda := i_X\omega\vert_\Sigma$ we compute
    \begin{align*}
        i_X\omega\vert_\Sigma \wedge \omega\vert_\Sigma^{n - 1} = \frac{1}{n} i_X\omega\vert_\Sigma^n.
    \end{align*}
    As $X$ is nowhere tangent to $\Sigma$, the form $i_X\omega\vert_\Sigma^n \in \Omega^{2n - 1}( \Sigma)$ is a volume form by \cite[Proposition~15.21]{lee:dt:2012}. Finally, using Cartan's magic formula \ref{lem:cartan}, we compute
    \begin{equation*}
        \ker \omega\vert_\Sigma \subseteq \ker L_X\omega\vert_\Sigma = \ker (i_Xd\omega + di_X\omega)\vert_\Sigma = \ker d\lambda.
    \end{equation*}
\end{solution}

\begin{example}[{Magnetic Torus, \cite[Section~6.1]{cieliebakfrauenfelderpaternain:mane:2010}}]
	\label{ex:magnetic_torus}
	Let $\mathbb{T}^n$ be the standard flat torus for $n \geq 2$ and let $J \colon \mathbb{R}^n \to \mathbb{R}^n$ be  an antisymmetric nonzero linear map. Define $\rho \in \Omega^2(\mathbb{T}^n)$ by setting $\rho(\cdot,\cdot) := \langle \cdot, J \cdot \rangle$ and denote by $\omega_\rho = dp \wedge dq + \pi^*\rho$ the magnetic symplectic form on $T^*\mathbb{T}^n \cong \mathbb{T}^n \times \mathbb{R}^n$. For an energy value $c \in \mathbb{R}$ set $\Sigma_c := H^{-1}(c)$ for the mechanical Hamiltonian function 
	\begin{equation*}
		H(q,p) := \frac{1}{2} \norm[0]{p}^2 \qquad \forall (q,p) \in \mathbb{T}^n \times \mathbb{R}^n.
	\end{equation*}
	Define $A := (J\vert_{\im J})^{-1}$ and $\alpha \in \Omega^1(\im J)$ by
	\begin{equation*}
		\alpha_x(v) := \frac{1}{2}\langle x,Av\rangle.
	\end{equation*}
	Then $\Sigma_c$ is a stable Hamiltonian manifold for every $c > 0$ by \cite[Proposition~6.3]{cieliebakfrauenfelderpaternain:mane:2010}. The stabilising form $\lambda$ on $\Sigma_c$ is given by
	\begin{equation}
		\label{eq:stabilising_form}
		\lambda = f^*(pdq) + (\pr_\parallel \circ \pr)^*\alpha,
	\end{equation}
	\noindent where 
	\begin{equation*}
		\pr_\perp \colon \mathbb{R}^n \to \ker J, \qquad \pr_\parallel \colon \mathbb{R}^n \to \im J, \qquad \pr \colon \mathbb{T}^n \times \mathbb{R}^n \to \mathbb{R}^n
	\end{equation*}
	\noindent denote the projections with respect to the orthogonal splitting
	\begin{equation*}
		\mathbb{R}^n = \ker J \oplus \im J,
	\end{equation*}
	\noindent and
	\begin{equation*}
		f \colon \mathbb{T}^n \times \mathbb{R}^n \to \mathbb{T}^n \times \mathbb{R}^n, \qquad f(q,p) := \del[1]{q,\pr_\perp(p)}.
	\end{equation*}
\end{example}

\begin{definition}[{Contact Manifold, \cite[Definition~2.5.1]{frauenfelderkoert:3bp:2018}}]
	A \bld{contact manifold}\index{Manifold!contact} is defined to be a stable Hamiltonian manifold $(\Sigma,d\lambda,\lambda)$. We simply write $(\Sigma,\lambda)$ for a contact manifold and call $\lambda$ a \bld{contact form}.
\end{definition}

\begin{example}[{Regular Energy Surface, \cite[Theorem~1.2.2]{abbashofer:cg:2019}}]
	Let $\Sigma$ be a compact regular energy surface in a mechanical Hamiltonian system $(T^*M,dp \wedge dq,H)$. Then there exists $\lambda \in \Omega^1(\Sigma)$ such that $d\lambda = \omega\vert_\Sigma$ and $(\Sigma,\lambda)$ is a contact manifold.
\end{example}

\begin{definition}[{Reeb Vector Field, \cite[p.~1773]{cieliebakfrauenfelderpaternain:mane:2010}}]
	\label{def:Reeb}
	Let $(\Sigma,\omega,\lambda)$ be a stable Hamiltonian manifold. The unique vector field $R \in \mathfrak{X}(\Sigma)$ implicitly defined by 
	\begin{equation*}
		i_R\omega = 0 \qquad \text{and} \qquad i_R \lambda = 1
	\end{equation*}
	\noindent is called the \bld{Reeb vector field}\index{Vector field!Reeb}.
\end{definition}

\begin{example}[Star-Shaped Hypersurfaces]
    \label{ex:star-shaped}
    Let $f \in C^\infty(\mathbb{S}^{2n - 1})$ be a positive function. Then the star-shaped hypersurface
    \begin{equation*}
        \Sigma_f = \{f(z)z : z \in \mathbb{S}^{2n - 1}\}
    \end{equation*}
    \noindent is a contact manifold with contact form $\lambda\vert_{\Sigma_f}$, where
    \begin{equation*}
        \lambda := \frac{1}{2}\sum_{j = 1}^n(y_j dx_j - x_j dy_j).
    \end{equation*}
    Indeed, by \cite[Lemma~12.2.2]{frauenfelderkoert:3bp:2018}, we have that
    \begin{equation*}
        X_{H_f}\vert_{\Sigma_f} \in \ker d\lambda\vert_{\Sigma_f} \qquad \text{and} \qquad \lambda(X_{H_f})\vert_{\Sigma_f} = 1
    \end{equation*}
    \noindent for the defining Hamiltonian function
    \begin{equation*}
        H_f \colon \mathbb{C}^n \setminus \{0\} \to \mathbb{R}, \qquad H_f(z) := \frac{\norm[0]{z}^2}{f(z/\norm[0]{z}^2)} - 1.
    \end{equation*}
    Hence it follows from Example \ref{ex:stable} that $(\Sigma_f,\lambda\vert_{\Sigma_f})$ is a contact manifold as the vector field
    \begin{equation*}
        X := \frac{1}{2}\sum_{j = 1}^n\del[3]{x_j\frac{\partial}{\partial x_j} + y_j \frac{\partial}{\partial y_j}} \in \mathfrak{X}(\mathbb{R}^{2n})
    \end{equation*}
    \noindent satisfies $i_X d\lambda = \lambda$ and is nowhere tangent to $\Sigma$ as
    \begin{equation*}
        dH_f(X)\vert_{\Sigma_f} = d\lambda(X,X_{H_f})\vert_{\Sigma_f} = \lambda(X_{H_f})\vert_{\Sigma_f} = 1.
    \end{equation*}
    Finally, we conclude that
    \begin{equation*}
        X_{H_f}\vert_{\Sigma_f} = R_f \in \mathfrak{X}(\Sigma_f)
    \end{equation*}
    \noindent is the Reeb vector field.
\end{example}

\begin{figure}[h!tb]
	\centering
	\includegraphics[width=.6\textwidth]{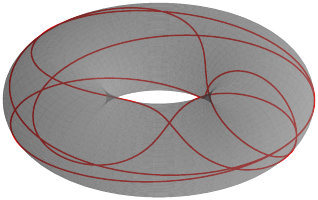}
	\caption{A contractible periodic Reeb orbit on a magnetic torus $\mathbb{T}^2$}
	\label{fig:torus}
\end{figure}

\begin{example}[{Magnetic Torus, \cite[Proposition~6.3]{cieliebakfrauenfelderpaternain:mane:2010}}]
	\label{ex:magnetic_flow}
	The flow $\phi_t$ of the magnetic Hamiltonian system in Example \ref{ex:magnetic_torus} is given by 
	\begin{equation*}
		\phi_t(q,p) =\del[3]{\int_0^t e^{sJ}p ds + q, e^{tJ}p},
	\end{equation*}
	\noindent as one can explicitly compute this flow using \eqref{eq:magnetic_Hamiltonian_vf}, and $(q,p) \in \Sigma_c$ gives rise to a contractible closed orbit of period $\tau$ if and only if
	\begin{equation*}
		\label{eq:magnetic_flow}
		\pr_\perp(p) = 0, \qquad e^{\tau J} \pr_\parallel(p) = \pr_\parallel(p), \qquad \text{and} \qquad \norm[1]{\pr_\parallel(p)}^2 = 2c.
	\end{equation*}
	It is illustrative to consider the special case $n = 2$ and $\sigma = dq_1 \wedge dq_2$. Then
	\begin{equation*}
		\lambda = -\frac{1}{2}(p_1dp_2 - p_2dp_1),
	\end{equation*}
	\noindent and the projection of a contractible periodic orbit to $\mathbb{T}^2$ is depicted in Figure \ref{fig:torus}.
\end{example}

\begin{exercise}
	Let $(\Sigma,\omega,\lambda)$ be a compact stable Hamiltonian manifold. 
    \begin{enumerate}[label=\alph*.]
        \item Show that the flow $\phi_t \colon \Sigma \times \mathbb{R} \to \Sigma$ of the Reeb vector field $R \in \mathfrak{X}(\Sigma)$ preserves $\lambda$, that is, it holds that
	       \begin{equation*}
		      \phi^*_t \lambda = \lambda \qquad \forall t \in \mathbb{R}.
	        \end{equation*}
         \item Explain, why we defined stable Hamiltonian manifolds as we did in Definition \ref{def:stable}.
         \item Conclude, that in the case where $\lambda$ is a contact form, the linearised flow $D\phi_t$ preserves the splitting
	\begin{equation*}
		T\Sigma = \ker \lambda \oplus \langle R \rangle.
	\end{equation*}
    \end{enumerate}
\end{exercise}

\begin{solution}
    \mbox{}
    \begin{enumerate}[label=\alph*.]
        \item Using Fisherman's formula \ref{lem:fisherman} and Cartan's magic formula \ref{lem:cartan} we compute
    \begin{equation*}
        \frac{d}{dt}\phi_t^*\lambda = \phi_t^* L_R \lambda = \phi_t^*(i_R d\lambda + di_R \lambda) = 0
    \end{equation*}
    \noindent as $R \in \ker \omega \subseteq \ker d\lambda$ and $i_R\lambda = 1$. Thus $\phi_t^*\lambda$ is constant and so $\phi_t^*\lambda = \phi_0^* \lambda = \lambda$.
        \item In the definition of the Reeb vector field we see that we need the condition of $\lambda$ being nowhere vanishing on the characteristic distribution $\ker \omega$ to have existence and uniqueness. The condition $\ker \omega \subseteq \ker d\lambda$ ensures that part a is true, which has important dynamical consequences.
        \item By part a we have that $\lambda = \phi^*_t \lambda = \lambda \circ D\phi_t$. As $D\phi_t$ is an isomorphism, we conclude $\ker \lambda = \ker \lambda \circ D\phi_t$. 
    \end{enumerate}
\end{solution}
\newpage
\section{The Limit Set of a Family of Periodic Orbits}
In this final section we answer the question
\begin{center}
    \emph{What is the dynamical meaning of stability?}
\end{center}
This question was partially answered in \cite[Lemma~2.5]{cieliebakfrauenfelderpaternain:mane:2010} and we will generalise \cite[Theorem~A]{belbrunofrauenfelderkoert:omega-limit_set:2020}, which excludes blue sky catastrophes in the stable case.

\subsection{Blue Sky Catastrophes}

\begin{definition}[Homotopy of Stable Energy Surfaces]
	\label{def:stable_homotopy}
	Let $(M,\omega)$ be a connected symplectic manifold. A \bld{homotopy of stable energy surfaces} is defined to be a time-dependent Hamiltonian function $H \in C^\infty(M \times [0,1])$ such that
	\begin{enumerate}
		\item $0$ is a regular value of $H_\sigma := H(\cdot,\sigma) \in C^\infty(M)$ for all $\sigma \in [0,1]$.
		\item $\Sigma_\sigma := H_\sigma^{-1}(0)$ is connected for all $\sigma \in [0,1]$.
		\item $H^{-1}(0)$ is compact.
		\item There exists a smooth time-dependent vector field $X$ on a neighbourhood of $H^{-1}(0)$ such that
			\begin{equation*}
				dH_\sigma(X_\sigma)\vert_{\Sigma_\sigma} > 0 \qquad \text{and} \qquad \ker(\omega\vert_{\Sigma_\sigma}) \subseteq \ker(L_{X_\sigma}\omega\vert_{\Sigma_\sigma}) 
			\end{equation*}
			\noindent for all $\sigma \in [0,1]$.
	\end{enumerate}
	We also write $(\Sigma_\sigma,\omega_\sigma,\lambda_\sigma)_{\sigma \in [0,1]}$ for a homotopy of stable regular energy surfaces, where
	\begin{equation*}
		\omega_\sigma := \omega\vert_{\Sigma_\sigma} \qquad \text{and} \qquad \lambda_\sigma := i_{X_\sigma}\omega\vert_{\Sigma_\sigma}
	\end{equation*}
	\noindent for all $\sigma \in [0,1]$.
\end{definition}

\begin{figure}[h!tb]
	\centering
	\includegraphics[width=.8\textwidth]{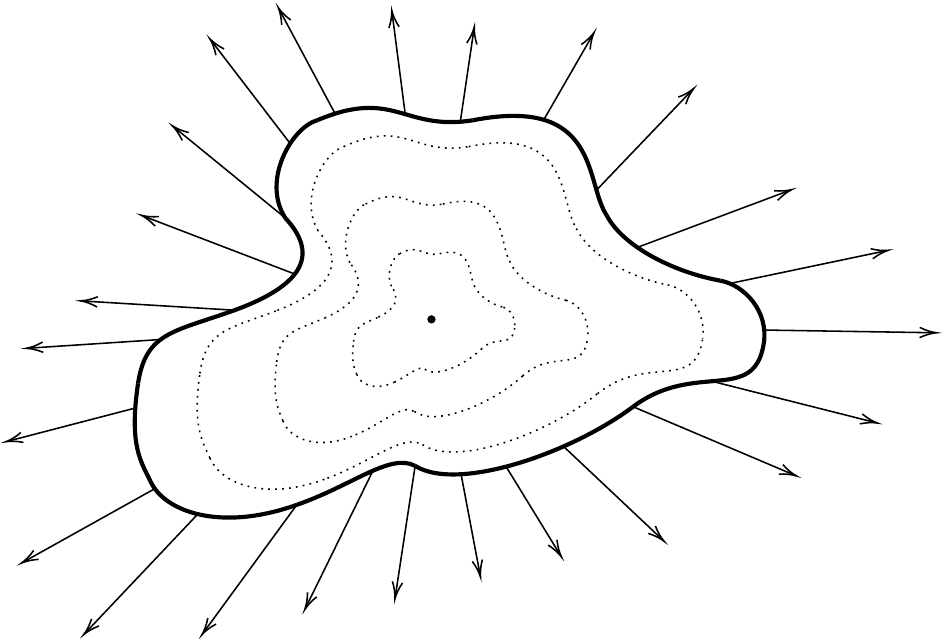}
	\caption{A stable homotopy of star-shaped hypersurfaces in $\mathbb{C}^n$}
	\label{fig:star-shaped}
\end{figure}

\begin{example}[Star-Shaped Hypersurfaces]
	Let $(f_\sigma)_{\sigma \in [0,1]}$ be a smooth family of functions in $C^\infty(\mathbb{S}^{2n - 1})$. Then the smooth family $(H_\sigma)_{\intcc[0]{0,1}} \subseteq C^\infty(\mathbb{C}^n \setminus \{0\})$
    \begin{equation*}
		H_\sigma(z) := \frac{\norm[0]{z}^2}{f_\sigma(z/\norm[0]{z}^2)} - 1, \qquad \sigma \in \intcc[0]{0,1}
    \end{equation*}
    \noindent defines a homotopy of stable energy surfaces. Indeed, by Example \ref{ex:star-shaped}, every hypersurface $(\Sigma_{f_\sigma},\lambda\vert_{\Sigma_{f_\sigma}})$ is a contact manifold.
\end{example}
    
\begin{example}[Magnetic Torus]
    Let $0 < a < b$ and consider the magnetic Hamiltonian system from Example \ref{ex:magnetic_torus}. Then the smooth family $H_\sigma := H - (a + \sigma (b - a))$ is a homotopy of stable energy surfaces.
\end{example}

Let $(\Sigma_\sigma,\omega_\sigma,\lambda_\sigma)_{\sigma \in \intcc[0]{0,1}}$ be a homotopy of stable energy surfaces. Consider a smooth family of parametrised periodic orbits $(\gamma_\sigma,\tau_\sigma)_{\sigma \in \intco[0]{0,\sigma_\infty}}$ for $0 < \sigma_\infty \leq 1$, that is, 
\begin{equation*}
	(\gamma_\sigma,\tau_\sigma) \in \mathscr{L}\Sigma_\sigma \times \intoo[0]{0,+\infty}
\end{equation*}
\noindent solves the equation
\begin{equation}
    \label{eq:parametrised}
	\dot{\gamma}_\sigma(t) = \tau_\sigma X_{H_\sigma}(\gamma_\sigma(t)) \qquad \forall t \in \mathbb{T},
\end{equation}
\noindent for all $\sigma \in \intco[0]{0,\sigma_\infty}$. 

\begin{example}[Star-Shaped Hypersurfaces]
    By \cite[Theorem~1.1]{abreu:symmetric:2022}, there does exist a closed Reeb orbit on any star-shaped hypersurface $\Sigma_f$, giving rise to a smooth family of parametrised periodic orbits. This does finally answer Exercise \ref{ex:orbits_star-shaped} part d. Another approach to this question is via Rabinowitz--Floer homology in \cite[Corollary~1.5]{cieliebakfrauenfelder:rfh:2009} .
\end{example}

\begin{example}[Magnetic Torus]
    By Example \ref{ex:magnetic_flow} there does exist a closed Reeb orbit on any hypersurface $\Sigma_c$ of the magnetic torus, giving rise to a smooth family of parametrised periodic orbits. More abstractly, the existence of such orbits is also given by Schlenk's Theorem \cite[Theorem~4.9]{cieliebakfrauenfelderpaternain:mane:2010}.
\end{example}

In order to systematically study the limit behaviour of the smooth family 
$(\gamma_\sigma,\tau_\sigma)_{\sigma \in \intco[0]{0,\sigma_\infty}}$ it is convenient to introduce the \bld{limit set} consisting of all 
\begin{equation*}
	(\gamma_{\sigma_\infty},\tau_{\sigma_\infty}) \in \mathscr{L}\Sigma_{\sigma_\infty} \times \intoo[0]{0,+\infty}
\end{equation*}
\noindent such that there exists a sequence $(\sigma_k) \subseteq \intco[0]{0,\sigma_\infty}$ with
\begin{equation*}
	\sigma_k \to \sigma_\infty \qquad \text{and} \qquad (\gamma_{\sigma_k},\tau_{\sigma_k}) \xrightarrow{C^\infty} (\gamma_{\sigma_\infty},\tau_{\sigma_\infty})
\end{equation*}
\noindent as $k \to \infty$. If $\tau_\sigma$ is uniformly bounded on $\intco[0]{0,\sigma_\infty}$, the limit set is not empty by Ascolis theorem and bootstrapping \eqref{eq:parametrised}. However, it could be that $\tau_\sigma \to +\infty$ as $\sigma \to \sigma_\infty$. This scenario is referred to as a \bld{blue sky catastrophe}, see \cite[p.~121]{frauenfelderkoert:3bp:2018}. In the case of a homotopy of stable  energy surfaces, no blue sky catastrophe does occur because of the following result.

\begin{theorem}
    \label{thm:main}
	Let $(\gamma_\sigma,\tau_\sigma)_{\sigma \in \intco[0]{0,\sigma_\infty}}$ be a smooth family of parametrised periodic orbits for some $0 < \sigma_\infty \leq 1$ on a homotopy of stable  energy surfaces $(\Sigma_\sigma,\omega_\sigma,\lambda_\sigma)_{\sigma \in [0,1]}$. Then there exists a constant $C > 0$ such that
	\begin{equation*}
		\frac{1}{C} \leq \tau_\sigma \leq C \qquad \forall \sigma \in \intco[0]{0,\sigma_\infty}.
	\end{equation*}
    In particular, the limit set of $(\gamma_\sigma,\tau_\sigma)_{\sigma \in \intco[0]{0,\sigma_\infty}}$ is not empty.
\end{theorem}

\begin{corollary}
    The limit set of a smooth family of parametrised periodic orbits on a homotopy of stable regular energy surfaces is nonempty, compact and connected. 
\end{corollary}

\begin{proof}
    That the limit set is nonempty is Theorem \ref{thm:main}. That it is compact and connected goes along the lines of the proof \cite[Theorem~5.2]{belbrunofrauenfelderkoert:omega-limit_set:2020}.
\end{proof}

\begin{definition}[Rabinowitz Action Functional]
    Let $(M,\lambda,H)$ be an exact Hamiltonian system. Then the corresponding \bld{Rabinowitz action functional} is
\begin{equation*}
    \mathscr{A}^H \colon \mathscr{L}M \times \mathbb{R} \to \mathbb{R}, \qquad \mathscr{A}^H(\gamma,\tau) := \int_0^1 \gamma^*\lambda - \tau\int_0^1 H(\gamma(t))dt.
\end{equation*}
\end{definition}

\begin{exercise}
    \label{ex:compare}
    Let $(M,\lambda,H)$ be an exact Hamiltonian system.
    \begin{enumerate}[label=\alph*.]
        \item What are the differences and similarities of the Hamiltonian action functional $\mathscr{A}_H$ and the Rabinowitz action functional $\mathscr{A}^H$? Do sketch the corresponding critical points!
        \item Why is the functional $\mathscr{A}^H$ better suited for the proof of Theorem \ref{thm:main} than $\mathscr{A}_H$? 
    \end{enumerate}
\end{exercise}

\begin{figure}[h!tb]
     \centering
     \begin{subfigure}[b]{0.45\textwidth}
         \centering
         \includegraphics[width=\textwidth]{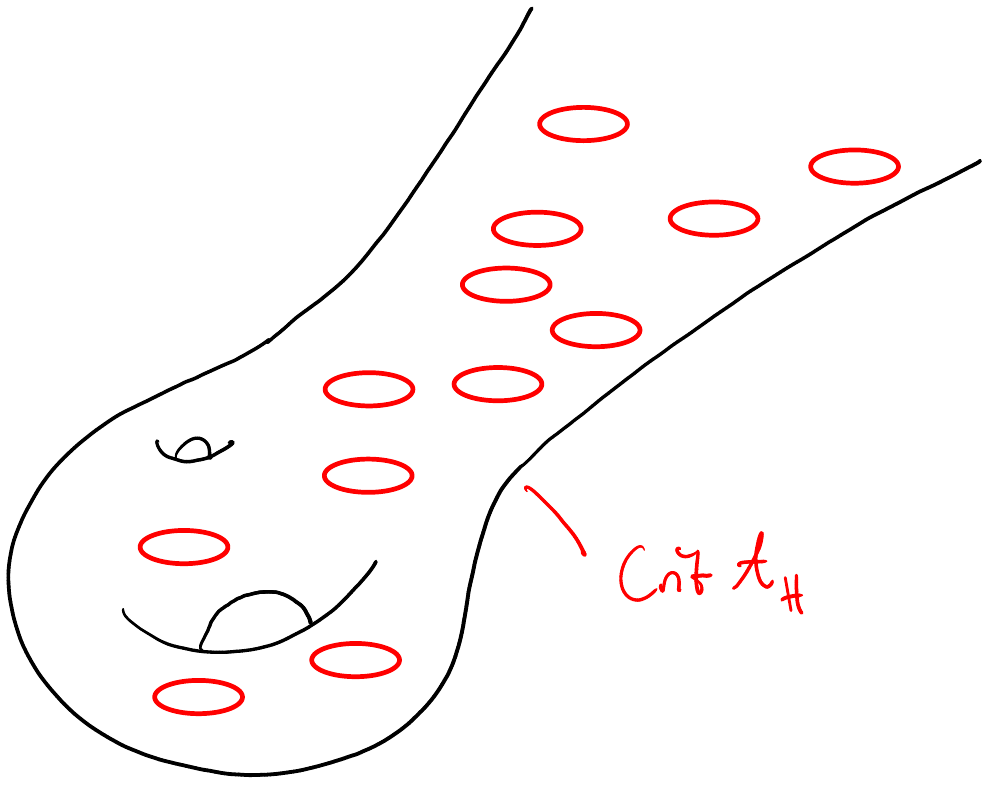}
     \end{subfigure}
     \hfill
     \begin{subfigure}[b]{0.45\textwidth}
         \centering
         \includegraphics[width=\textwidth]{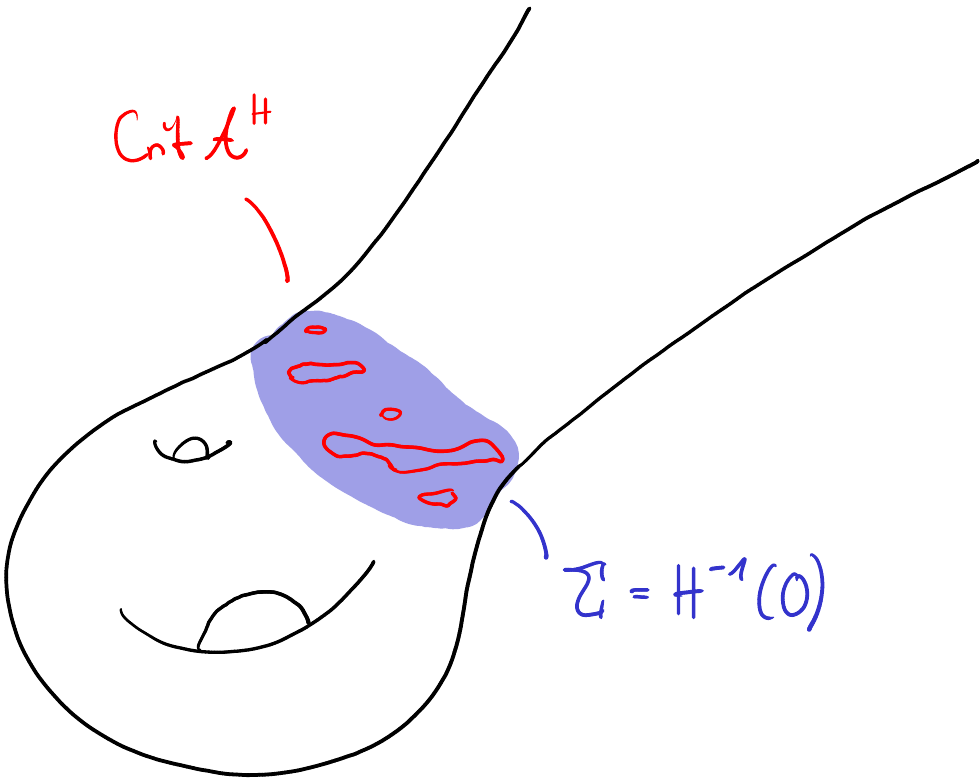}
     \end{subfigure}
\end{figure}

\begin{solution}
    \mbox{}
    \begin{enumerate}[label=\alph*.]
        \item Both functionals are defined on infinite-dimensional spaces and differ only by the Lagrange multiplicator $\tau$ in the Rabinowitz action functional $\mathscr{A}^H$. Thus the critical points of $\mathscr{A}^H$ are $\tau$-periodic loops of constant energy of the Hamiltonian vector field $X_H$, whereas the critical points of the Hamiltonian action functional $\mathscr{A}_H$ are $1$-periodic loops of arbitrary energy. Hence the functional $\mathscr{A}_H$ is better suited for problems involving fixed period but arbitrary energy, and the Rabinowitz action functional is better suited for problems involving fixed energy but arbitrary period. To summarise, we have that
        \begin{equation*}
            \gamma \in \Crit(\mathscr{A}_H) \qquad \Leftrightarrow \qquad \dot{\gamma}(t) = X_H(\gamma(t)) \> \forall t \in \mathbb{T},
        \end{equation*}
        \noindent and
        \begin{equation*}
            (\gamma,\tau) \in \Crit(\mathscr{A}^H) \qquad \Leftrightarrow \qquad \begin{cases}
                \gamma \in \mathscr{L}\Sigma,\\
                \dot{\gamma}(t) = X_H(\gamma(t)) \> \forall t \in \mathbb{T}.
            \end{cases} 
        \end{equation*}
        The proof for the description of the critical points of the Rabinowitz action function is similar to the proof of Lemma \ref{lem:hamiltonian_equations} and makes use of preservation of energy \ref{cor:preservation_of_energy}.
        
        \item The Rabinowitz action functional is better suited as the statement of Theorem \ref{thm:main} involves a fixed energy but arbitrary period problem. However, it is not at all clear, how one should start a proof of this result and why it is an advantage to use the notion of a Rabinowitz action functional.
    \end{enumerate}   
\end{solution}

\begin{exercise}
    Read the following proof of Theorem \ref{thm:main}. What does simplify if we assume that the homotopy is a homotopy of contact hypersurfaces? You can also compare the proofs given in \cite[Section~7.6]{frauenfelderkoert:3bp:2018} and \cite{belbrunofrauenfelderkoert:omega-limit_set:2020}.
\end{exercise}

\begin{proof}[Proof of Theorem \ref{thm:main}]
    We proceed in two steps.
    
	\emph{Step 1: The family $(\gamma_\sigma,\tau_\sigma)_{\sigma \in \intco[0]{0,\sigma_\infty}}$ gives rise to a smooth family of critical points of suitable Rabinowitz action functionals.} Let $(\varepsilon,\psi_\sigma)_{\sigma \in \intcc[0]{0,1}}$ be a smooth family of stable tubular neighbourhoods of $(\Sigma_\sigma,\omega_\sigma,\lambda_\sigma)$ as in \cite[Proposition~2.6~(a)]{cieliebakfrauenfelderpaternain:mane:2010}. Fix mollifications $f \in C^\infty(\mathbb{R})$ and $h \in C^\infty(\mathbb{R})$ such that
	\begin{itemize}
		\item $f(r) = r + 1$ for $r \in \intcc[1]{-\frac{\varepsilon}{2},\frac{\varepsilon}{2}}$.
		\item $\supp f \subseteq \intoo[0]{-\varepsilon,\varepsilon}$.
		\item $h$ satisfies
			\begin{equation*}
				h(r) = \begin{cases}
					r & r \in \intcc[1]{-\frac{\varepsilon}{3},\frac{\varepsilon}{3}}\\
					\frac{\varepsilon}{3} & r \in \intco[1]{\frac{\varepsilon}{3},+\infty}\\
					-\frac{\varepsilon}{3} & r \in \intoc[1]{-\infty,-\frac{\varepsilon}{3}}
				\end{cases}.
			\end{equation*}
	\end{itemize}
	See Figure \ref{fig:mollifications}.

	\begin{figure}[h!tb]
	\begin{subfigure}[c]{0.5\textwidth}
		\centering
		\includegraphics[width=.8\textwidth]{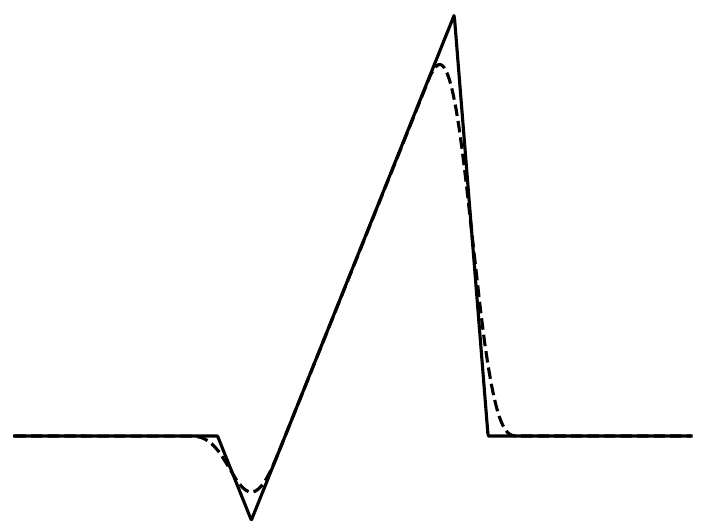}
		\subcaption{Possible choice of extension function $f$ and its mollification}
		\label{fig:f}
	\end{subfigure}
	~
	\begin{subfigure}[c]{0.5\textwidth}
		\centering
		\includegraphics[width=.8\textwidth]{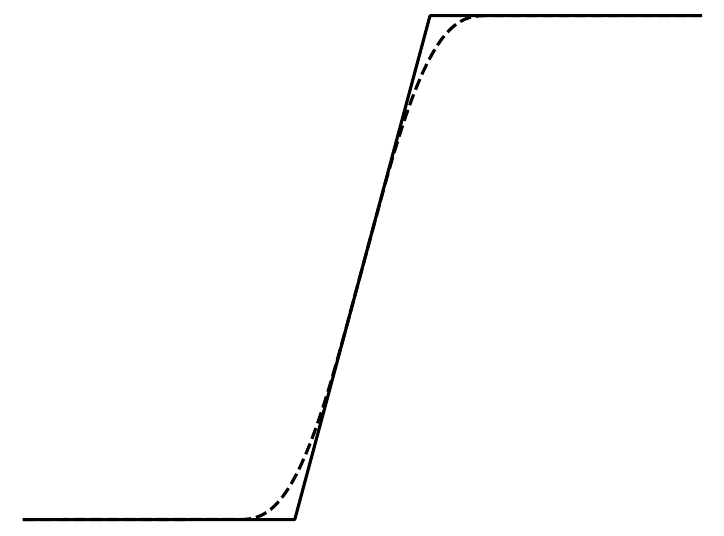}
		\subcaption{Mollification of the piecewise-linear function $h$}
		\label{fig:h}
	\end{subfigure}
	\caption{Mollifications}
	\label{fig:mollifications}
\end{figure}

Define extensions 
\begin{equation*}
	\bar{\lambda}_\sigma(p) := \begin{cases}
		f(r)\lambda_\sigma(x) & p = \psi_\sigma(r,x),\\
		0 & p \notin U_\sigma,
	\end{cases}
\end{equation*}
\noindent and modifications
\begin{equation*}
	\bar{H}_\sigma(p) := \begin{cases}
		h(r) & p = \psi_\sigma(r,x),\\
		\frac{\varepsilon}{3} & p \in U_\sigma^c \cap H^{-1}_\sigma(0,+\infty),\\
		-\frac{\varepsilon}{3} & p \in U_\sigma^c \cap H^{-1}_\sigma(-\infty,0),\\
	\end{cases}
\end{equation*}
\noindent for $U_\sigma := \psi_\sigma(\intoo[0]{-\varepsilon,\varepsilon} \times \Sigma_\sigma)$. Note that the definition of $\bar{H} \in C^\infty(M \times \intcc[0]{0,1})$ makes sense, as by assumption, every $\Sigma_\sigma$ is a separating hypersurface. For every $\sigma \in \intcc[0]{0,1}$, consider the Rabinowitz action functional
\begin{equation*}
	\mathscr{A}_\sigma \colon \mathscr{L} M \times \intoo[0]{0,+\infty} \to \mathbb{R}
\end{equation*}
\noindent defined by
\begin{equation*}
	\mathscr{A}_\sigma(\gamma,\tau) := \int_0^1 \gamma^*\bar{\lambda}_\sigma - \tau\int_0^1 \bar{H}_\sigma \circ \gamma.
\end{equation*}
We claim that
\begin{equation}
    \label{eq:critical_points_inclusion}
	d\mathscr{A}_\sigma\vert_{(\gamma,\tau)}(X,\eta) = \bar{m}\del[1]{\grad \mathscr{A}_\sigma\vert_{(\gamma,\tau)},(X,\eta)}
\end{equation}
\noindent holds for all $(X,\eta) \in T_\gamma\mathscr{L}M \times \mathbb{R}$ and $(\gamma,\tau) \in \mathscr{L} M \times \mathbb{R}$, where
\begin{equation*}
	\bar{m}((X,\eta),(Y,\xi)) := \int_0^1 d\bar{\lambda}_\sigma(JX,Y) + \eta\xi,
\end{equation*}
\noindent for all $(Y,\xi) \in \mathscr{L} M \times \mathbb{R}$ with $J$ being an $\omega$-compatible almost complex structure on $(M,\omega)$, and
\begin{equation*}
	\grad \mathscr{A}_\sigma\vert_{(\gamma,\tau)} := \begin{pmatrix}
		J(\dot{\gamma} - \tau X_{\bar{H}_\sigma} \circ \gamma)\\
		\displaystyle-\int_0^1 \bar{H}_\sigma \circ \gamma
	\end{pmatrix}.
\end{equation*}
We compute
\allowdisplaybreaks
\begin{align*}
	i_{h'(r)R_\sigma}\psi^*_\sigma \omega &= i_{h'(r)R_\sigma}(\omega_\sigma + d(r\lambda_\sigma))\\
	&= i_{h'(r)R_\sigma}(\omega_\sigma + dr \wedge \lambda_\sigma + rd\lambda_\sigma)\\
	&=h'(r)i_{R_\sigma}\omega_\sigma + h'(r)i_{R_\sigma}dr \lambda_\sigma - h'(r)i_{R_\sigma}\lambda_\sigma dr + rh'(r)i_{R_\sigma}d\lambda_\sigma\\
	&= -h'(r)dr\\
	&= -d\bar{H}_\sigma,
\end{align*}
\noindent as $R_\sigma$ belongs to $\ker \omega_\sigma \subseteq \ker d\lambda_\sigma$ by stability, and
\allowdisplaybreaks
\begin{align*}
	i_{X_{\bar{H}_\sigma}}d\bar{\lambda}_\sigma &= i_{X_{\bar{H}_\sigma}}(f'(r)dr \wedge \lambda_\sigma + f(r) d\lambda_\sigma)\\
	&= f'(r)i_{X_{\bar{H}_\sigma}}dr \lambda_\sigma - f'(r)i_{X_{\bar{H}_\sigma}}\lambda_\sigma dr + f(r)i_{X_{\bar{H}_\sigma}}d\lambda_\sigma\\
	&= -f'(r)h'(r)dr\\
	&= -h'(r)dr\\
	&= -d\bar{H}_\sigma.
\end{align*}
Hence
\begin{align*}
	d\mathscr{A}_\sigma\vert_{(\gamma,\tau)}(X,\eta) &= \int_0^1 d\bar{\lambda}_\sigma(X,\dot{\gamma}) - \tau \int_0^1 d\bar{H}_\sigma(X) - \eta\int_0^1 \bar{H}_\sigma \circ \gamma\\
	&= \int_0^1 d\bar{\lambda}_\sigma(X,\dot{\gamma} - \tau X_{\bar{H}_\sigma} \circ \gamma)- \eta\int_0^1 \bar{H}_\sigma \circ \gamma\\
	&= \bar{m}\del[1]{\grad \mathscr{A}_\sigma\vert_{(\gamma,\tau)},(X,\eta)}.
\end{align*}

\emph{Step 2: The periods  $\tau_\sigma$ are uniformly bounded from below and above.} As in \cite[Theorem~7.6.1]{frauenfelderkoert:3bp:2018}, we have the period-action equality 
\begin{align*}
	\mathscr{A}_\sigma(\gamma_\sigma,\tau_\sigma) &= \int_0^1 \gamma_\sigma^* \bar{\lambda}_\sigma - \tau_\sigma \int_0^1 \bar{H}_\sigma \circ \gamma_\sigma\\
	&= \int_0^1 \bar{\lambda}_\sigma(\dot{\gamma}_\sigma)\\
	&= \int_0^1 \lambda_\sigma(\dot{\gamma}_\sigma)\\
	&= \tau_\sigma\int_0^1 \lambda_\sigma(R_\sigma \circ \gamma_\sigma)\\
	&= \tau_\sigma.
\end{align*}
Using Step 1 we compute
\begin{align*}
	\partial_\sigma\tau_\sigma &= \partial_\sigma(\mathscr{A}_\sigma(\gamma_\sigma,\tau_\sigma))\\
	&= (\partial_\sigma \mathscr{A}_\sigma)(\gamma_\sigma,\tau_\sigma) + d\mathscr{A}_\sigma\vert_{(\gamma_\sigma,\tau_\sigma)}(\partial_\sigma(\gamma_\sigma,\tau_\sigma))\\
	&= (\partial_\sigma \mathscr{A}_\sigma)(\gamma_\sigma,\tau_\sigma) + \bar{m}\del[1]{\grad \mathscr{A}_\sigma\vert_{(\gamma_\sigma,\tau_\sigma)},(\partial_\sigma\gamma_\sigma,\partial_\sigma\tau_\sigma)}\\
	&=  (\partial_\sigma \mathscr{A}_\sigma)(\gamma_\sigma,\tau_\sigma)\\
	&= \int_0^1 (\partial_\sigma \bar{\lambda}_\sigma)(\dot{\gamma}_\sigma) - \tau_\sigma\int_0^1 (\partial_\sigma \bar{H}_\sigma) \circ \gamma_\sigma\\
	&= \tau_\sigma \int_0^1 (\partial_\sigma\bar{\lambda}_\sigma)(R_\sigma \circ \gamma_\sigma) - \tau_\sigma \int_0^1 (\partial_\sigma\bar{H}_\sigma) \circ \gamma_\sigma
\end{align*}
\noindent because $\grad \mathscr{A}_\sigma\vert_{(\gamma_\sigma,\tau_\sigma)} = 0$ for all $\sigma \in \intco[0]{0,\sigma_\infty}$. As $H^{-1}(0)$ is compact by assumption, there exists $\kappa > 0$ such that
\begin{equation*}
	\abs[0]{(\partial_\sigma \bar{\lambda}_\sigma)(R_\sigma(x))} \leq \frac{\kappa}{2} \qquad \text{and} \qquad \abs[0]{\partial_\sigma \bar{H}_\sigma(x)} \leq \frac{\kappa}{2}.
\end{equation*}
\noindent for all $(x,\sigma) \in H^{-1}(0)$. Consequently, we have that $\abs[0]{\partial_\sigma \tau_\sigma} \leq \kappa\tau_\sigma$ and so
\begin{equation*}
	0 < \tau_{\sigma_0}e^{-\kappa(\sigma_1 - \sigma_0)} \leq \tau_{\sigma_1} \leq \tau_{\sigma_0}e^{\kappa(\sigma_1 - \sigma_0)} < +\infty
\end{equation*}
\noindent for all $0 \leq \sigma_0 < \sigma_1 < \sigma_\infty$. Thus
\begin{equation*}
	\frac{1}{C} \leq \tau_\sigma \leq C \qquad C := \tau_0 e^{\kappa\sigma_\infty}
\end{equation*}
\noindent for all $\sigma \in \intco[0]{0,\sigma_\infty}$.
\end{proof}

\begin{question}
    In the setting of Theorem \ref{thm:main} and its proof, define the set
    \begin{equation*}
        \mathscr{P}(\Sigma_\sigma) := \{(\gamma,\tau) \in \mathscr{L}\Sigma_\sigma \times \intoo[0]{0,+\infty} : \dot{\gamma}(t) = \tau R_\sigma(\gamma(t))\>\forall t \in \mathbb{T}\}
    \end{equation*}
    \noindent of \bld{parametrised periodic Reeb orbits on $\Sigma_\sigma$}. Which statements are true?
    \begin{itemize}
        \item[$\square$] $\Crit(\mathscr{A}_\sigma) \subseteq \mathscr{P}(\Sigma_\sigma)$ 
        \item[$\square$] $\mathscr{P}(\Sigma_\sigma) \subseteq \Crit(\mathscr{A}_\sigma)$ 
        \item[$\square$] $\mathscr{P}(\Sigma_\sigma) = \Crit(\mathscr{A}_\sigma)$ 
        \item[$\square$] There is no relation between $\mathscr{P}(\Sigma_\sigma)$ and $\Crit(\mathscr{A}_\sigma)$ as $X_{\bar{H}_\sigma}\vert_{\Sigma_\sigma} \neq R_\sigma$.
    \end{itemize}
\end{question}

\begin{solution}
    The only correct item is the second one which follows from \eqref{eq:critical_points_inclusion}. The reason why the third item is not true in general, is that $d\bar{\lambda}_\sigma$ is degenerate in most cases. For example, the stabilising form could be closed.
\end{solution}

\begin{exercise}
    Answer the question \emph{what is the dynamical meaning of stability?}
\end{exercise}

\begin{solution}
    The dynamical meaning of stability is the fact that the periods of parametrised Reeb orbits on a homotopy of stable energy surfaces remain stable. More precisely, the period is uniformly bounded from above and below by Theorem \ref{thm:main}. The stability property is crucially used in the proof of Theorem \ref{thm:main}.
\end{solution}

\begin{remark}
    In the restricted three-body problem families of periodic orbits are known to exist \cite{belbruno:lunar:2019}.
\end{remark}

\subsection{Equivariant Rabinowitz Action Functionals} 
Consider the $\mathbb{Z}_m$-action on the odd-dimensional sphere $\mathbb{S}^{2n - 1}$ induced by the rotation
\begin{equation*}
    \varphi \colon \mathbb{C}^n \to \mathbb{C}^n, \qquad \varphi(z_1,\dots,z_n) := \del[1]{e^{2\pi i k_1/m}z_1,\dots, e^{2\pi i k_n/m}z_n},
\end{equation*}
\noindent where $m \geq 1$ is an integer and $k_1,\dots,k_n \in \mathbb{Z}$ are coprime to $m$. The resulting smooth quotient manifold $\mathbb{S}^{2n - 1}/\mathbb{Z}_m$ is called a \emph{lens space}. The contact form $\lambda\vert_{\mathbb{S}^{2n - 1}}$ from Example \ref{ex:star-shaped} is invariant under this $\mathbb{Z}_m$-action and thus descends to the quotient to a contact form on the lens space. The existence of closed Reeb orbits on lens spaces is important in the study of celestial mechanics. Indeed, by \cite[Corollary~5.7.5]{frauenfelderkoert:3bp:2018}, the Moser regularised energy hypersurface near the earth or the moon of the planar circular restricted three-body problem for energy values below the first critical value is diffeomorphic to the real projective space $\mathbb{RP}^3$. 

\begin{definition}[$G$-Invariant Loop Space]
    Let $G$ be a symmetry group of a Hamiltonian system $(M,\omega,H)$. Define the \bld{$G$-invariant loop space of $M$} by
    \begin{equation*}
        \mathscr{L}_G M := \{\gamma \in C^\infty(\mathbb{T},M) : \psi_g(\gamma(\mathbb{T})) = \gamma(\mathbb{T}) \> \forall g \in G\}.
    \end{equation*}
\end{definition}

\begin{example}[Twisted Loop Space]
    Let $\varphi \in \Diff(M)$ be of finite order $m := \ord(\varphi)$ such that $\varphi^* \lambda = \lambda$ and $H \circ \varphi = H$ holds for an exact Hamiltonian system $(M,\lambda,H)$. We consider the induced $\mathbb{Z}_m$-action on $M$ given by
    \begin{equation*}
        \mathbb{Z}_m \times M \to M, \qquad [k]\cdot x = \varphi^k(x).
    \end{equation*}
    Then $\mathbb{Z}_m$ is a symmetry group of $(M,d\lambda,H)$, and
    \begin{equation*}
        \mathscr{L}_\varphi M = \{\gamma \in C^\infty(\mathbb{R},M) : \gamma(t + 1) = \varphi(\gamma(t))\>\forall t \in \mathbb{R}\} \subseteq \mathscr{L}_{\mathbb{Z}_m} M
    \end{equation*}
    \noindent is the \bld{twisted loop space of $M$ and $ \varphi$}.
\end{example}

\begin{exercise}
    Does the definition of a $G$-invariant loop space make sense for a continuous symmetry group $G$?
\end{exercise}

\begin{solution}
    For a general continuous symmetry group the definition makes no sense. Consider for example the Kepler problem \ref{ex:kepler}. There we have $\mathscr{L}_{\operatorname{SO}(n)}T^*(\mathbb{R}^n \setminus\{0\}) = \varnothing$.
\end{solution}

\begin{definition}[{Twisted Rabinowitz Action Functional, \cite[Definition~2.6]{baehni:rfh:2023}}]
    Let $(M,\lambda,H)$ be an exact Hamiltonian system. For a diffeomorphism $\varphi \in \Diff(M)$ with $\varphi^*\lambda = \lambda$ and $H \circ \varphi = H$, we define the \bld{twisted Rabinowitz action functional}
    \begin{equation*}
        \mathscr{A}_\varphi^H \colon \mathscr{L}_\varphi M \times \mathbb{R} \to \mathbb{R}, \qquad \mathscr{A}_\varphi^H(\gamma,\tau) := \int_0^1 \gamma^*\lambda - \tau \int_0^1 H(\gamma(t))dt.
    \end{equation*}
\end{definition}

Using the twisted Rabinowitz action functional one can prove the following result.

\begin{theorem}[{\cite[Theorem~1.2]{baehni:rfh:2023}}]
    \label{thm:noncontractible}
	Let $\Sigma \subseteq \mathbb{C}^n$, $n \geq 2$, be a compact and connected star-shaped hypersurface invariant under the rotation	
	\begin{equation*}
		\varphi \colon \mathbb{C}^n \to \mathbb{C}^n, \quad \varphi(z_1,\dots,z_n) := \del[1]{e^{2\pi i k_1/m}z_1, \dots, e^{2\pi i k_n/m}z_n}
	\end{equation*}
	\noindent for some even integer $m \geq 2$ and $k_1,\dots,k_n \in \mathbb{Z}$ coprime to $m$. Then the quotient $\Sigma/\mathbb{Z}_m$ admits a noncontractible periodic Reeb orbit generating the fundamental group $\pi_1(\Sigma/\mathbb{Z}_m) \cong \mathbb{Z}_m$.
\end{theorem}

\begin{remark}
    In \cite{abreu:symmetric:2022} the restriction of $m$ being even was removed and an additional multiplicity result was given \cite[Theorem~1.2]{abreu:symmetric:2022}. As the construction works for any $\mathbb{Z}_m$-action, one might consider a full $\mathbb{S}^1$-action. The definition of the $\mathbb{S}^1$-equivariant Rabinowitz action functional is given in \cite{frauenfelderschlenk:equivariant:2016}. Yet another similar variant was used in \cite[Section~4.1]{albers:gysin:2023}.
\end{remark}

So what about continuous symmetries? One could make an attempt and adapt the functional for defining moment Floer homology in \cite[Section~4.2]{frauenfelder:moment:2004}.

\begin{definition}[Moment Rabinowitz Action Functional]
    For a symmetry group $G$ of an exact Hamiltonian system $(M,\lambda,H)$ we define the \bld{moment Rabinowitz action functional}
    \begin{equation*}
        \mathscr{A}^{H,\mu} \colon C^\infty(\mathbb{T}, M \times \mathfrak{g}) \times \mathbb{R} \to \mathbb{R}
    \end{equation*}
    \noindent by
    \begin{equation*}
        \mathscr{A}^{H,\mu}(\gamma,\xi,\tau) := \int_0^1 \gamma^*\lambda - \tau \int_0^1 H(\gamma(t))dt - \int_0^1 \mu(\xi(t))(\gamma(t))dt.
    \end{equation*}
\end{definition}

\begin{exercise}
    Is it worth studying this functional?
\end{exercise}

\begin{solution}
    If the symmetry is discrete, then the functional $\mathscr{A}^{H,\mu}$ reduces to the standard Rabinowitz action functional $\mathscr{A}^H$. If $G$ is not discrete, the critical points of the moment Rabinowitz action functional are given by
    \begin{equation*}
        (\gamma,\xi,\tau) \in \Crit(\mathscr{A}^{H,\mu}) \qquad \Leftrightarrow \qquad \begin{cases}
            \gamma \in \mathscr{L}\Sigma,\\
            \dot{\gamma}(t) = \tau X_H(\gamma(t)) - \widehat{\xi}(t) \> \forall t \in \mathbb{T},\\
            \mu^*(\gamma(t)) = 0 \> \forall t \in \mathbb{T},
        \end{cases}
    \end{equation*}
    \noindent where
    \begin{equation*}
        \mu^* \colon M \to \mathfrak{g}^*, \qquad \mu^*(x)(\xi) := \mu(\xi)(x)
    \end{equation*}
    \noindent denotes the comoment map. Thus if we pass to the quotient, we consider periodic orbits on
    \begin{equation*}
        \del[1]{H^{-1}(0) \cap (\mu^*)^{-1}(0)}/G \cong H^{-1}(0)/G \cap (\mu^*)^{-1}(0)/G.
    \end{equation*}
    Under additional assumptions on the moment map $\mu$, the quotient $(\mu^*)^{-1}(0)/G$ is the well-known Marsden--Weinstein quotient, a compact symplectic manifold. This symplectic form cannot be exact because of compactness (see Question \ref{que:symplectic}). Thus it seems difficult to investigate this scenario with the pictures from Exercise \ref{ex:compare} in mind.
\end{solution}

\section*{Acknowledgements}
To Colin, Neil and Jil. In memory of Will J.
Merry. A brilliant teacher and a guiding light. Without him I would have never met Urs Frauenfelder, Kai Cieliebak and Felix Schlenk.
\newpage
\addcontentsline{toc}{section}{References}
\printbibliography

\end{document}